     \documentclass[reqno]{amsart}
\usepackage{amsmath}
\usepackage{amsfonts}
\usepackage{amstext}
\usepackage{amsbsy}
\usepackage{amsopn}
\usepackage{amsxtra}
\usepackage{upref}
\usepackage{amsthm}
\usepackage{amsmath}
\usepackage{amssymb}

\usepackage[textures,bookmarks=false]{hyperref}
\usepackage{mathrsfs}
\usepackage{enumerate}

\theoremstyle{plain}
\newtheorem{theorem}{Theorem}
\newtheorem{maintheorem}{Theorem}

\newtheorem{maincorollary}[maintheorem]{Corollary}
\newtheorem{lemma}{Lemma}
\newtheorem{proposition}{Proposition}
\newtheorem{remark}{Remark}

\def\bd{\partial}
\def\es{\emptyset}
\def\sm{\setminus}
\def\st{such that }
\def\R{\mathbb{R}}
\def\N{\mathbb{N}}
\def\crit{\mbox{\rm Crit}}
\def\F{\mathcal{F}}

\def\eps{\varepsilon}
\def\phi{\varphi}
\def\le{\leqslant}
\def\ge{\geqslant}
\def\c{{\bf C}}
\def\T{\mathcal{T}}

\def\htop{h_{top}}
\def\ie{{\em i.e.,~}}
\def\M{\mathcal{M}}
\def\uG{\underline{LG}}
\def\oG{\overline{LG}}
\def\P{\mathcal{P}}
\def\diam{\text{diam}}
\def\D{\mathcal{D}}
\def\freq{\text{freq}}

\def\level{\text{lev}}
\def\ae{a.e.\ }

\def\CCE{C \negthinspace_{ \scriptscriptstyle C\negmedspace E}}
\def\betaCE{\beta \negthinspace _{\scriptscriptstyle C\negmedspace E}}
\def\CP{C\negthinspace_{\scriptscriptstyle P}>0}
\def\betaP{\beta\negthinspace_{\scriptscriptstyle P}}

\def\DS{\mathcal{DS}}
\def\K{\mathcal{K}}

\begin{document}

\title{Multifractal analysis for multimodal maps}
\date{}
\author{Mike Todd}
\address{Centro de Matem\'atica da Universidade do Porto, Rua do Campo Alegre 687, 4169-007 Porto, Portugal \footnote{
{\bf Current address:}\\
Department of Mathematics and Statistics\\
Boston University\\
111 Cummington Street\\
Boston, MA 02215\\
USA }\\ }
\email{mtodd@math.bu.edu}
\urladdr{http://math.bu.edu/people/mtodd/}


\thanks{This work was supported by FCT grant SFRH/BPD/26521/2006 and also by FCT through CMUP}

\begin{abstract}
Given a multimodal interval map $f:I \to I$ and a H\"older potential $\phi:I \to \R$, we study the dimension spectrum for equilibrium states of $\phi$.  The main tool here is inducing schemes, used to overcome the presence of critical points.  The key issue is to show that enough points are `seen' by a class of inducing schemes.  We also compute the Lyapunov spectrum.  We obtain the strongest results when $f$ is a Collet-Eckmann map, but our analysis also holds for maps satisfying much weaker growth conditions along critical orbits.
\end{abstract}

\maketitle

\section{Introduction}

Given a metric space $X$ and a probability measure $\mu$ on $X$, the \emph{pointwise dimension} of $\mu$ at $x\in X$ is defined as
$$d_\mu(x):=\lim_{r\to 0^+}\frac{\log\mu(B_r(x))}{\log r}$$
if the limit exists, where $B_r(x)$ is a ball of radius $r$ around $x$.  This tells us how concentrated a measure is around a point $x$\iffalse.
For example, if $\mu$ is an atomic measure supported at a point $x$, then $d_\mu(x)=0$.  However, if $y\neq x$ then $d_\mu(y)=\infty$.\fi;
the more concentrated, the lower the value of $d_\mu(x)$.
For an endomorphism $f:X\to X$, we will study the pointwise dimension of $f$-invariant measures $\mu$.  In particular we will be interested in equilibrium states $\mu_\phi$ for $\phi:X\to \R$ in a certain class of potentials (see below for definitions).

For any $A\subset X$, we let $\dim_H(A)$ denote the Hausdorff dimension of $A$. We let $$\K_\phi(\alpha):=\left\{x:\lim_{r\to 0^+}\frac{\log\mu_{\phi}(B_r(x))}{\log r} =\alpha\right\}, \quad \DS_\phi(\alpha):=\dim_H(\K_\phi(\alpha)),$$ and
$$\K_\phi':=\left\{x:\lim_{r\to 0^+}\frac{\log\mu_{\phi}(B_r(x))}{\log r} \text{ does not exist}\right\}.$$  Then we can make a \emph{multifractal decomposition}:
$$X=\K'_\phi\cup\left(\cup_{\alpha\in \R}\K_\phi(\alpha)\right).$$
The function $\DS_\phi$ is known as the \emph{dimension spectrum} of $\mu_\phi$.  The study of this function fits into the more general theory of thermodynamic formalism which also gives us information on the statistical properties of the system such as return time statistics, large deviations and decay of correlations.

These ideas are generally well understood in the case of uniformly hyperbolic systems, see \cite{Pesbook}.  The dimension spectrum can be described in terms of the pressure function, which we define below.  A common way to prove this in uniformly hyperbolic cases is to code the system using a finite Markov shift, and then exploit the well developed theory of thermodynamic formalism and dimension spectra for Markov shifts, see for example \cite{PesWeiss}.  For non-uniformly hyperbolic dynamical systems this approach can be more complicated since we often need to code by countable Markov shifts.  As has been shown by Sarig \cite{Sarphase, Sarcts}, Iommi \cite{Iom, Iomrenyi} and Pesin and Zhang \cite{PeZphase} among others, in going from finite to countable Markov shifts, more exotic behaviour, including `phase transitions', appears.

The coding used in non-uniformly hyperbolic cases usually arises from an `inducing scheme': that is, for some part of the phase space, iterates of the original map are taken, and the resulting `induced map' is considered.  The induced maps are Markov, and so the theory of countable Markov shifts as in \cite{HMU, Iom} can be used.  In some cases the induced map can be a first return map to an interval, but this is not always so.

There has been a lot of success with the inducing approach in the case of Manneville-Pomeau maps.  These are interval maps which are expanding everywhere, except at a parabolic fixed point.  The presence of the parabolic point leads to phase transitions as mentioned above.  Multifractal analysis, of the dimension spectrum and the Lyapunov spectrum (see below), of these examples has been carried out by Pollicott and Weiss \cite{PolWei}, Nakaishi \cite{Nak} and Gelfert and Rams \cite{GelRams}.  In the first two of these papers, inducing schemes were used (in the third one, the fact that the original system is Markov is used extensively).  The inducing schemes used are first return maps to a certain natural domain.  The points of the original phase space which the inducing schemes do not `see' is negligible, consisting only of the (countable) set preimages of the parabolic point.  We also mention a closely related theory for certain Kleinian groups by Kesseb\"ohmer and Stratmann \cite{KesStr}.

In the case of multimodal maps with critical points, if the critical orbits are dense then there is no way that useful inducing schemes can be first return maps to intervals.  Moreover, the set of points which the inducing schemes do not `see' can, in principle, be rather large.  %(As explained in Section~\ref{sec:what is seen}, this set will at %least include the (countable) set of points which eventually map to %a critical point.)
In these cases the thermodynamic formalism has a lot of exotic behaviour: phase transitions brought about due to some polynomial growth condition were discussed by Bruin and Keller in \cite{BrKel} and shown in more detail by Bruin and Todd \cite{BTeqnat}.  Multiple phase transitions, which are due to renormalisations rather than any growth behaviour, were proved by Dobbs \cite{Dobbsphase}.

In this paper we develop a multifractal theory for maps with critical points by defining inducing schemes which provide us with sufficient information on the dimension spectrum.  The main idea is that points with large enough pointwise Lyapunov exponent must be `seen' by certain inducing schemes constructed in \cite{BTeqnat}.  These inducing schemes are produced via the Markov extension known as the Hofbauer extension, also known as the Hofbauer tower.  This structure was developed by Hofbauer and Keller, see for example \cite{Htop, Hpwise, Kellift}.  Their principle applications were for interval maps.  The theory for higher dimensional cases was further developed by Buzzi \cite{Bumulti}.  Once we have produced these inducing schemes, we can use the theory of multifractal analysis developed by Iommi in \cite{Iom} for the countable Markov shift case.  Note that points with pointwise Lyapunov exponent zero cannot be `seen' by measures which are compatible to an inducing scheme, so if we are to use measures and inducing schemes to study the dimension spectrum, the inducing methods presented here may well be optimal.

There is a further property which useful inducing schemes must have: not only must they see sufficiently many points, but also they must be well understood from the perspective of the thermodynamic formalism.  Specifically, given a potential $\psi$, we need its induced version on the inducing scheme to fit into the framework of Sarig \cite{SaBIP}.  In \cite{PeSe, BTeqgen, BTeqnat} this was essentially translated into having `good tail behaviour' of the equilibrium states for the induced potentials.

Our main theorem states that, as in the expanding case, for a large class of multimodal maps, the multifractal spectrum can be expressed in terms of the Legendre transform of the pressure function for important sets of parameters $\alpha$.  The Collet-Eckmann case is closest to the expanding case, and here we indeed get exactly the same kind of graph for $\alpha\mapsto \DS_\phi(\alpha)$ as in the expanding case for the values of $\alpha$ we consider. In the non-Collet Eckmann case, we expect the graph of $\DS_\phi$ to be qualitatively different from the expanding case, as shown for the related Lyapunov spectrum in \cite{Nak} and \cite{GelRams}.  We note that singular behaviour of the Lyapunov spectrum was also observed by Bohr and Rand \cite{BoRa} for the special case of the quadratic Chebyshev polynomial.

The results presented here can be seen as an extension of some of the ideas in \cite{Hofdim}, in which the full analysis of the dimension spectrum was only done for uniformly expanding interval maps.  See also \cite{Yurimulti} for maps with weaker expansion properties.  Moreover, Hofbauer, Raith and Steinberger \cite{HoRaSt} proved the equality of various thermodynamic quantities for non-uniformly expanding interval maps, using `essential multifractal dimensions'.  However, the full analysis in the non-uniformly expanding case, including the expression of the dimension spectrum in terms of some Legendre transform, was left open.

\subsection{Key definitions and main results}

Given a dynamical system $f:X\to X$, we let
$$\M=\M(f):=\{\text{$f$-invariant probability measures on } X\}$$  and  $$\M_{erg}=\M_{erg}(f):=\{\mu\in \M: \mu \text{ is } \text{ergodic}\}.$$
For a potential $\phi:X\to \R$, the \emph{pressure} is defined as
$$P(\phi):=\sup_{\mu\in \M}\left\{h_\mu+\int\phi~d\mu:-\int\phi~d\mu<\infty \right\}$$ where $h_\mu$ denotes the metric entropy with respect to $\mu$.  Note that by the ergodic decomposition, we can just take the above supremum over $\M_{erg}$.  We let $\htop(f)$ denote the topological entropy of $f$, which is equal to $P(0)$, see \cite{Kbook}.  A measure $\mu$ which `achieves the pressure', \ie $h_\mu+\int\phi~d\mu=P(\phi)$, is called an \emph{equilibrium state}.

\iffalse
We will be interested in $C^2$ multimodal interval maps $f:I\to I$.  Let $\crit=\crit (f)$ denote the set of critical points of $f$. We say that $c\in \crit$ is a
\emph{non-flat} critical point of $f$ if there exists a
diffeomorphism $g_c:\R \to \R$ with $g_c(0)=0$ and $1<\ell_c<\infty$
\st for $x$ close to $c$, $f(x)=f(c)\pm|g_c(x-c)|^{\ell_c}$.  The
value of $\ell_c$ is known as the \emph{critical order} of $c$.  We define $\ell_{max}(f):=\max\{\ell_c:c\in {\rm Crit}(f)\}$.  Our class of maps $\F$ will have all critical points non-flat, as well as some other properties we describe in more detail below.
\fi

Let $\mathcal F$ be the collection of $C^3$ multimodal interval maps $f:I \to I$ where $I=[0,1]$, satisfying:
\iffalse
\begin{itemize}
\item[(a)] The critical set $\crit = \crit(f)$ consists of finitely many critical points $c$ with critical orders $1 < \ell_c < \infty$, i.e., $f(x) = f(c) + (g(x-c))^{\ell_c}$ for some diffeomorphisms $g:\R \to \R$ with $g(0) = 0$ and $x$ close to $c$.
\item[(b)] $f$ has negative Schwarzian derivative, i.e., $1/\sqrt{|Df|}$ is convex.
\item[(c)] The non-wandering set $\Omega$ (the set of points $x\in I$ \st for arbitrarily small neighbourhoods $U$ of $x$ there exists $n=n(U)\ge 1$ \st $f^n(U)\cap U\neq \es$) consists of a single interval.
\item[(d)] $f^n(\crit)  \cap f^m(\crit)=\es$ for $m \neq n$.
\end{itemize}

\fi
\newcounter{Lcount}
\begin{list}{\alph{Lcount})}
{\usecounter{Lcount} \itemsep 1.0mm \topsep 0.0mm \leftmargin=7mm}

\item the critical set $\crit = \crit(f)$ consists of finitely many critical point $c$ with critical order $1 < \ell_c < \infty$, i.e., $f(x) = f(c) + (g(x-c))^{\ell_c}$ for some diffeomorphisms $g:\R \to \R$ with $g(0) = 0$ and $x$ close to $c$;
%\item $f$ has negative Schwarzian derivative, i.e., $1/\sqrt{|Df|}$ is convex;
\item $f$ has no parabolic cycles;
\item $f$ is topologically transitive on $I$;
\item $f^n(\crit)  \cap f^m(\crit)=\es$ for $m \neq n$.
\end{list}

\begin{remark}
Conditions c) and d) are for ease of exposition, but not crucial.  In particular, Condition c) excludes that $f$ is renormalisable.  For multimodal maps satisfying a) and b), the set $\Omega$ consists of finitely many components $\Omega_k$, on each of which $f$ is topologically transitive, see \cite[Section III.4]{MSbook}.  In the case where there is more than one transitive component in $\Omega$, for example the renormalisable case, the analysis presented here can be applied to any one of the transitive components consisting of intervals.   We also note that in this case $\Omega$ contains a (hyperbolic) Cantor set outside components of $\Omega$ which consist of intervals.  The work of Dobbs \cite{Dobbsphase} shows that renormalisable maps these hyperbolic Cantor sets can give rise to singular behaviour in the thermodynamic formalism (phase transitions in the pressure function $t\mapsto P(t\phi)$) not accounted for by the behaviour of critical points themselves.  For these components we could apply a version of the usual hyperbolic theory to study the dimension spectra.

We include condition b) in order to apply the distortion theorem,
\cite[Theorem C]{SVarg}.  Alternatively, we could assume negative Schwarzian derivative, since this added to the transitivity assumption implies that there are no parabolic points.

Condition d) rules out one critical point mapping onto another. Alternatively, it would be possible to consider these critical points as a `block', but to simplify the exposition, we will not do that here.  Condition d) also rules out critical points being preperiodic.

\end{remark}

We define the \emph{lower/upper pointwise Lyapunov exponent} as  $$\underline\lambda_f(x):=\liminf_{n\to\infty} \frac{1}{n} \sum_{j=0}^{n-1} \log|Df(f^j(x))|, \text{ and } \overline\lambda_f(x):=\limsup_{n\to\infty} \frac{1}{n} \sum_{j=0}^{n-1} \log|Df(f^j(x))|$$ respectively.  If $\underline\lambda_f(x)=\overline\lambda_f(x)$, then we write this as $\lambda_f(x)$.  For a measure $\mu\in \M_{erg}$, we let $$\lambda_f(\mu):=\int\log|Df|~d\mu$$ denote the Lyapunov exponent of the measure.  Since our definition of $\F$ will exclude the presence of attracting cycles, \cite{Prz} implies that $\lambda_f(\mu)\ge 0$ for all $f\in\F$ and $\mu\in \M$.

For $\lambda\ge 0$, we denote the `good Lyapunov exponent' sets by $$\uG_\lambda:=\{x:\underline\lambda_f(x)>\lambda\} \text{ and }
\oG_\lambda:=\{x:\overline\lambda_f(x)>\lambda\}.$$
We define $$\tilde\K_\phi(\alpha):=\K_\phi(\alpha)\cap \oG_0 \text{ and } \widetilde{\DS}_\phi(\alpha):=\dim_H(\tilde\K_\phi(\alpha)).$$

As well as assuming that our maps $f$ are in $\F$, we will also sometimes impose certain growth conditions on $f$:

\begin{list}{$\bullet$}{\itemsep 1.0mm \topsep 0.0mm } \setlength{\itemindent=-5mm}%\leftmargin=7mm}
\item An exponential growth condition (Collet-Eckmann): there exist $\CCE, \betaCE>0$,
\begin{align}  |Df^n(f(c))|
\ge C\negthinspace_{\scriptscriptstyle  C\negmedspace E}e^{\beta\negthinspace _{\scriptscriptstyle  C\negmedspace E} n} \hbox{ for all } c\in \crit \hbox{ and } n\in \N.
\label{eq:CE}
\end{align}

\item A polynomial growth condition:
There exist $\CP>0$ and $\betaP>2\ell_{max}(f)$ so that
\begin{align}
|Df^n(f(c))|
\ge C\negthinspace_{\scriptscriptstyle P}n^{\beta\negthinspace_{\scriptscriptstyle P}} \hbox{ for all } c\in \crit \hbox{ and } n\in \N.
\label{eq:poly}
\end{align}

\item A simple growth condition:
\begin{align}
|Df^n(f(c))| \to \infty \hbox{ for all } c\in \crit.
\label{eq:grow}
\end{align}
\end{list}

In all of these cases, \cite{BRSS} implies that there is a unique absolutely continuous invariant probability measure (acip).  This measure has positive entropy by \cite[Exercise V.1.4]{MSbook} and \cite[Proposition 7]{SVarg}.

We will consider potentials $-t\log|Df|$ and also $\epsilon$-H\"older potentials $\phi:I \to \R$ satisfying
\begin{align}\sup\phi-\inf\phi<\htop(f).\label{eq:range}\end{align}
Without loss of generality, we will also assume that $P(\phi)=0$.  Note that our results do not depend crucially on $\epsilon\in (0,1]$, so we will ignore the precise value of $\epsilon$ from here on.

\begin{remark}
We would like to emphasise that \eqref{eq:range} may not be easy to remove as an assumption on our class of H\"older potentials if all the results we present here are to go through.  For example, in the setting of Manneville-Pomeau maps, in \cite[Section 6]{BTeqgen} it was shown that for any $\eps>0$, there exists a H\"older potential $\phi$ with $\sup\phi-\inf\phi=\htop(f)+\eps$ and for which the equilibrium state is a Dirac measure on the fixed point (which is not seen by any inducing scheme).
\end{remark}

We briefly sketch some properties of these maps and potentials. For details, see Propositions~\ref{prop:eqgen} and \ref{prop:eqnat}.  As we will see below, we are interested in potentials of the form $-t\log|Df|+\gamma\phi$.
By \cite{BTeqnat} if $f$ satisfies \eqref{eq:CE} then there exist $t_1<1<t_2$  such that for each $t\in (t_1,t_2)$ there is an equilibrium state $\mu_{-t\log|Df|}$ for $-t\log|Df|$.  If $f$ only satisfies \eqref{eq:poly} then we take $t_2=1$.  Combining \cite{BTeqnat} and \cite{BTeqgen}, for H\"older potentials $\phi$ we have equilibrium states $\mu_{-t\log|Df|+\gamma\phi}$ for $-t\log|Df|+\gamma\phi$ if $t$ is close to 1 and $\gamma$ is close to 0.
\iffalse
Keller shows that for a piecewise continuous map $f:I \to I$ and $\phi:I \to \R$ satisfying \eqref{eq:range}, there is an equilibrium state $\mu_\phi$.  Note that by \cite{BTeqgen}, strengthening the conditions on $f$ allows us to get equilibrium states for more exotic potentials, see appendix.
\fi
Also, by \cite{BTeqgen}, if \eqref{eq:grow} holds and $\phi$ is a H\"older potential satisfying \eqref{eq:range}, then
there are equilibrium states $\mu_{-t\log|Df|+\gamma\phi}$ for $-t\log|Df|+\gamma\phi$ if $t$ is close to 0 and $\gamma$ is close to 1.  These equilibrium states are unique.  As explained in the appendix, \eqref{eq:grow} is assumed in \cite{BTeqgen} in order to ensure that the induced versions of $\phi$ are sufficiently regular, so if this regularity can be shown another way, for example in the simple case that $\phi$ is a constant everywhere, this condition can be omitted.

We define the auxiliary function $T_\phi(q)$ to be so that
\begin{equation} P(\psi_q)=0, \text{ where } \psi_q:=-T_\phi(q)\log|Df|+q\phi.\label{eq:Tphi}\end{equation}
The map $q\mapsto T_\phi(q)$ is convex and if $P(\phi)=0$ then  $T_\phi(1)=0$. By Ledrappier \cite[Theorem 3]{Ledrap}, if there is an acip then it is an equilibrium state for $x\mapsto -\log|Df(x)|$ and so $T_\phi(0)=1$.  It may be the case that for some values of $q$, there is no such number.  For example, let $f\in \F$ be a unimodal map not satisfying \eqref{eq:CE}.  Then as in \cite{NoSa}, $P(-t\log|Df|)=0$ for all $t\ge 1$.  If we set $\phi$ to be the constant potential, then $P(\phi)=0$ implies $\phi\equiv -\htop(f)$ since then $P(\phi)=P(0)-\htop(f)=0$.  For such $\phi$ and for $q<0$, then $T_\phi(q)$ must be undefined.

For $h$ a convex function, we say that $(h,g)$ form a \emph{Fenchel pair} if
$$g(p)=\sup_x\{px-h(x)\}.$$
In this case, $g$ is known as the \emph{Legendre-Fenchel} transform of $h$.  If $h$ is convex and $C^1$ then the function $g$ is called the \emph{Legendre transform} of $h$ and
$$g(\alpha)=h(q)+q\alpha \text{ were } q \text{ is such that } \alpha=-Dh(q).$$

If $f\in \F$ satisfies \eqref{eq:grow} then \cite{BRSS} guarantees the existence and uniqueness of an acip $\mu_{-\log|Df|}$ and we let $$\alpha_{ac}:=\frac{-\int\phi~d\mu_{-\log|Df|}}{\lambda_f(\mu_{-\log|Df|})}.$$

\begin{maintheorem}  Suppose that $f\in \F$ is a map satisfying \eqref{eq:grow} and $\phi:I \to \R$ is a H\"older potential satisfying \eqref{eq:range}, and with $P(\phi)=0$.  If the equilibrium state $\mu_\phi$ is not equal to the acip then there exist open sets $U,V \subset \R$ so that $T_\phi$ is differentiable on $V$ and for $\alpha\in U$, the dimension spectrum $\alpha\mapsto \widetilde{\DS}_\phi(\alpha)$ is the Legendre transform of $q\mapsto T_\phi(q)$.
Moreover,
\begin{itemize}
\item[(a)]  $U$
contains a neighbourhood of  $\dim_H(\mu_\phi)$, and $\widetilde{\DS}_\phi(\dim_H(\mu_\phi))=\dim_H(\mu_\phi)$;

\item[(b)] if $f$ satisfies \eqref{eq:poly}, then $U$
contains both a neighbourhood of $\dim_H(\mu_\phi)$,
and a one-sided neighbourhood of $\alpha_{ac}$, where $\widetilde{\DS}_\phi(\alpha_{ac})=1$;

\item[(c)] if $f$ satisfies \eqref{eq:CE}, then $U$
contains both a neighbourhood of  $\dim_H(\mu_\phi)$ and of $\alpha_{ac}$.
\end{itemize}
Furthermore, for all $\alpha\in U$ there is a unique equilibrium state $\mu_{\psi_q}$ for the potential $\psi_q$ so that $\mu_{\psi_q}(\tilde\K_\alpha)=1$, where $\alpha=-DT_\phi(q)$.  This measure has full dimension on $\tilde\K_\alpha$, \ie $\dim_H(\mu_{\psi_q})=\dim_H(\tilde\K_\alpha)$.
\label{thm:main spectrum}
\end{maintheorem}

Note that by Hofbauer and Raith \cite{HofRaidim}, $\dim_H(\mu_\phi)=\frac{h_{\mu_{\phi}}}{\lambda_f(\mu_{\phi})}$, and as shown by Ledrappier \cite[Theorem 3]{Ledrap},
$\dim_H(\mu_{-\log|Df|})=\frac{h_{\mu_{-\log|Df|}}}{\lambda_f(\mu_{-\log|Df|})}=1$.

In Section~\ref{sec:LE spec} we consider the situation where $\phi$ is the constant potential, which we recall that since $P(\phi)=0$, must be of the form $\phi\equiv -\htop(f)$.  In that setting, as noted above $T_\phi$ is not defined for $q<0$ when $f$ is unimodal and does not satisfy \eqref{eq:CE}.  Therefore, in that case we would expect $\widetilde{\DS}_\phi$ to behave differently to the expanding case for $\alpha >\alpha_{ac}$.  This is why we only deal with a one-sided neighbourhood of $\alpha_{ac}$ in (b). See also Remark~\ref{rmk:phase trans} for more information on this.

If, contrary to the assumptions of Theorem~\ref{thm:main spectrum}, $\mu_\phi=\mu_{-\log|Df|}$ then $\widetilde{\DS}_\phi(\alpha)$ is zero for every $\alpha\in \R$, except at $\alpha=\dim_H(\mu_\phi)$, where it takes the value 1.  As in Remark~\ref{rmk:preper} below, for multimodal maps $f$ and $\phi$ a constant potential, this only occurs when $f$ has preperiodic critical points, for example when $f$ is the quadratic Chebyshev polynomial.  In view of Liv\v{s}ic theory for non-uniformly hyperbolic dynamical systems, in particular the results in \cite[Section 5]{BrHoNic}, we expect this to continue to hold for more general H\"older potentials.

According to \cite{BrvS}, if \eqref{eq:CE} holds then there exists $\lambda>0$ so that the nonwandering set $\Omega$ is contained in $\oG_\lambda\cup\left(\cup_{n\ge 0}f^{-n}(\crit)\right)$.  Therefore we have the following corollary. Note that here the neighbourhood $U$ is as in case (c) of Theorem~\ref{thm:main spectrum}.

\begin{maincorollary}
Suppose that $f\in \F$ satisfies the Collet-Eckmann condition \eqref{eq:CE} and $\phi:I \to \R$ is a H\"older potential satisfying \eqref{eq:range} and with $P(\phi)=0$. If the equilibrium state $\mu_\phi$ is not equal to the acip then
there exist open sets $U,V \subset \R$ so that $T_\phi$ is differentiable on $V$, $U$ contains $\dim_H(\mu_\phi)$ %$\frac{h_{\mu_{\phi}}}{\lambda_f(\mu_{\phi})}$
and 1, and so that for $\alpha\in U$ the dimension spectrum $\DS_\phi(\alpha)$ is the Legendre transform of $T_\phi$.
\label{cor:CE spectrum}
\end{maincorollary}

In fact, to ensure that $\widetilde{\DS}_\phi(\alpha)=\DS_\phi(\alpha)$ it is enough to show that `enough points iterate into a finite set of levels of the Hofbauer extension infinitely often'.  As in \cite{Kellift}, one way of guaranteeing this is to show that a large proportion of the sets we are interested in `go to large scale' infinitely often.  Graczyk and Smirnov \cite{GrSm}
%Proposition 7.3
showed that for rational maps of the complex plane satisfying a summability condition, this is true.  Restricting their result to real polynomials, we have the following Corollary, which we explain in more detail in Section~\ref{sec:GrSm}.

\begin{maincorollary}
Suppose that $f\in \F$ extends to a polynomial on $\mathbb C$ with no parabolic points, all critical points in $I$, and satisfying \eqref{eq:poly}.  Moreover, suppose that $\phi:I \to \R$ is a H\"older potential satisfying \eqref{eq:range} and $P(\phi)=0$. If the equilibrium state $\mu_\phi$ is not equal to the acip then
there exist sets $U, V\subset \R$ such that $U$ contains a one-sided neighbourhood of $\alpha_{ac}$, $T_\phi$ is differentiable on $V$, and for $\alpha\in U$ the dimension spectrum $\DS_\phi(\alpha)$ is the Legendre transform of $T_\phi$.  Moreover, if $\dim_H(\mu_\phi)>\frac{\ell_{max}(f)}{\beta_P-1}$ then the same is true for any $\alpha$ in a neighbourhood of $\dim_H(\mu_\phi)$. %$\dim_H(\mu_\phi)=\frac{h_{\mu_{\phi}}}{\lambda_f(\mu_{\phi})}$.
\label{cor:poly spectrum}
\end{maincorollary}

Barreira and Schmeling \cite{BarSch} showed that in many situations the set $\K_\phi'$ has full Hausdorff dimension.  As the following proposition states, this is also the case in our setting.  The proof follows almost immediately from \cite{BarSch}, but we give some details in Section~\ref{sec:main thm}.

\begin{proposition}
Suppose that $f\in \F$ satisfies \eqref{eq:grow} and $\phi:I \to \R$ is a H\"older potential satisfying \eqref{eq:range} and with $P(\phi)=0$.  Then $\dim_H(\K_\phi')=1$.
\label{prop:nonreg}
\end{proposition}

Theorem~\ref{thm:main spectrum} also allows us to compute the Lyapunov spectrum.  The results in this case are in Section~\ref{sec:LE spec}.

For ease of exposition, in most of this paper the potential $\phi$ is assumed to be H\"older.  In this case existence of an equilibrium state $\mu_\phi$ was proved by Keller \cite{kellholder}.  However, as we show in the appendix, all the results here hold for a class of potentials $(SVI)$ considered in \cite{BTeqgen}.
Therefore, as an auxiliary result, we prove the existence of conformal measures $m_\phi$ for potentials $\phi\in SVI$.  Moreover, for the corresponding equilibrium states $\mu_\phi$, the density $\frac{d\mu_\phi}{dm_\phi}$ is uniformly bounded away from 0 and $\infty$.  This is used here in order to compare $d_{\mu_\Phi}(x)$ and $d_{\mu_\phi}(x)$, where $\mu_\Phi$ is the equilibrium state for an inducing scheme $(X,F)$, with induced potential $\Phi:X \to \R$ (see below for more details).  The equality of $d_{\mu_\Phi}(x)$ and $d_{\mu_\phi}(x)$ for $x\in X$ is not immediate in either the case $\phi$ is H\"older or the case $\phi$ satisfies $SVI$.
This is in contrast to the situation where the inducing schemes are simply first return maps, in which case $\mu_\Phi$ is simply a rescaling of the original measure $\mu_\phi$ and hence $d_{\mu_\Phi}(x)=d_{\mu_\phi}(x)$.  However, we will prove that for the inducing schemes used here, this rescaling property is still true of the conformal measures $m_\phi$ and $m_\Phi$, which then allows us to compare $d_{\mu_\Phi}(x)$ and $d_{\mu_\phi}(x)$.  It is interesting to note that the proof of existence of a conformal measure also goes through for potentials of the form $x\mapsto -t\log|Df(x)|$.
\iffalse
 since the measure for the inducing scheme $\mu_\Phi$ is not, as would be the case if the inducing schemes were simply first return maps, simply a rescaling of the original measure $\mu_\phi$.  On the other hand, we show that this rescaling property is true of the conformal measures, which then allows us to compare $d_{\mu_\Phi}(x)$ and $d_{\mu_\phi}(x)$.
\fi

\textbf{Note added in proof:}
After this work was completed, it was communicated to me that in a work in progress, J.\ Rivera-Letelier and W.\ Shen have proved that in fact for any $f\in \F$,
$$\dim_H\{x:\overline\lambda_f(x)=0\}=0.$$
This strengthens the results presented in this paper, allowing us to replace $\widetilde{\DS}_\phi(\alpha)$ with ${\DS}_\phi(\alpha)$ throughout.

\emph{Acknowledgements:}
 \iffalse I would like to thank G.\ Iommi, H.\ Bruin, T.\ Jordan and N.\ Dobbs for useful conversations and comments on earlier versions of this paper.\fi
I would like to thank G.\ Iommi, H.\ Bruin, T.\ Jordan and N.\ Dobbs for useful comments on earlier versions of this paper.  I would also like to thank them and D. Rand for fruitful conversations.

\section{The maps, the measures and the inducing schemes}
\label{sec:maps meas ind}

\iffalse
Let $f:I \to I$ be a $C^2$ multimodal map of the unit interval $I$.
Throughout $\F$ will be the collection of $C^2$
interval maps which have negative Schwarzian (that is, $1/{\sqrt{|Df|}}$ is convex away from critical points) and all critical points non-flat, and not points of inflection.
We will assume for ease of exposition that for all $n\neq m$, $f^n(\crit)\cap f^m(\crit)=\es$.  Note that the situation is simpler if critical points have finite orbits, and that if one critical orbit maps to another, it is possible to consider these critical points together as a `block', but to simplify the exposition we will not do that here.
We will also assume for simplicity that maps $f\in \F$ are non-renormalisable, see \cite{MSbook}, and do not have any attracting periodic points.  This implies that the \emph{non-wandering set} $\Omega$: the set of points $x\in I$ \st for arbitrarily small neighbourhoods $U$ of $x$ there exists $n=n(U)\ge 1$ \st $f^n(U)\cap U\neq \es$, is a finite union $\cup_k\Omega_k$ such that each $\Omega_k$ is a finite union of intervals such that $f:\Omega_k \to \Omega_k$ is topologically transitive.  For ease of exposition, we will assume that $\Omega$ has only one component for all maps in $\F$.

\fi

Let $(X,f)$ be a dynamical system and $\phi:X\to [-\infty,\infty]$ be a potential.  For use later, we let
$$S_n\phi(x):=\phi(x)+\cdots+\phi\circ f^{n-1}(x).$$
We say that a measure $m$, is \emph{conformal} for $(X,f,\phi)$ if $m(X)=1$, and for any Borel set $A$ so that $f:A \to f(A)$ is a bijection, $$m(f(A))=\int_Ae^{-\phi}~dm$$  (or equivalently, $dm(f(x))=e^{-\phi(x)} dm(x)$).

\subsection{Hofbauer extensions}
\label{subsec:Hofbauer}
We next define the Hofbauer extension, sometimes also known as a Hofbauer tower.  The setup we present here can be applied to general dynamical systems, since it only uses the structure of dynamically defined cylinders.  An alternative way of thinking of the Hofbauer extension specifically for the case of multimodal interval maps, which explicitly makes use of the critical set, is presented in \cite{BrBr}.

We first consider the dynamically defined cylinders.  We let $\P_0:=I$ and $\P_n$ denote the collection of maximal intervals $\c_n$ so that $f^n:\c_n\to f^n(\c_n)$ is a homeomorphism.  We let $\c_n[x]$ denote the member of $\P_n$ containing $x$.  If $x\in \cup_{n\ge 0}f^{-n}(\crit)$ there may be more than one such interval, but this ambiguity will not cause us any problems here.

The \emph{Hofbauer extension} is defined as $$\hat
I:=\bigsqcup_{k\ge 0}\bigsqcup_{\c_{k}\in \P_{k}}
f^k(\c_{k})/\sim$$ where $f^k(\c_{k})\sim
f^{k'}(\c_{k'})$ as components of the disjoint union $\hat I$ if $f^k(\c_{k})= f^{k'}(\c_{k'})$ as subsets in $I$.  Let
$\D$ be the collection of domains of $\hat I$ and $\pi:\hat
I \to I$ be the natural inclusion map.  A point $\hat x\in \hat I$ can
be represented by $(x,D)$ where $\hat x\in D$ for $D\in \D$ and
$x=\pi(\hat x)$.  Given $\hat x\in \hat I$, we can denote the domain $D\in \D$ it belongs to by $D_{\hat x}$.

The map $\hat f:\hat I \to \hat I$ is defined by
$$\hat f(\hat x) = \hat f(x,D) = (f(x), D')$$
if there are cylinder sets $\c_k \supset \c_{k+1}$ \st $x \in
f^k(\c_{k+1}) \subset f^k(\c_{k}) = D$ and $D' = f^{k+1}
(\c_{k+1})$.
In this case, we write $D \to D'$, giving $(\D, \to)$ the
structure of a directed graph.  Therefore, the map $\pi$
acts as a semiconjugacy between $\hat f$ and $f$: $$\pi\circ \hat
f=f\circ \pi.$$ % It is easy to check that there is a
%one-to-one correspondence between cylinder sets $\c_k \in \P_k$
%and $k$-paths $D_0 \to \dots \to D_{k}$ starting at the base
We denote the `base' of $\hat I$, the copy of $I$ in $\hat I$ by  $D_0$.  For $D\in \D$, we define $\level(D)$ to be the length of the shortest path $D_0 \to \dots \to D$ starting at the base $D_0$.  For each $R \in \N$, let $\hat I_R$ be the compact
part of the Hofbauer extension defined by the disjoint union
$$
\hat I_R := \sqcup \{ D \in \D : \level(D) \le R \}.$$

For maps in $\F$, we can say more about the graph structure of $(\D, \to)$ since Lemma 1 of \cite{BTeqnat} implies that if $f\in \F$ then there is a closed primitive subgraph $\D_{\T}$ of $\D$.  That is, for any $D,D' \in\D_{\T}$ there is a path $D\to \cdots \to D'$; and for any $D\in \D_{\T}$, if there is a path $D\to D'$ then $D'\in \D_{\T}$ too.  We can denote the disjoint union of these domains by $\hat I_{\T}$.  The same lemma says that if $f\in \F$ then $\pi(\hat I_{\T})=\Omega$ and $\hat f$ is transitive on $\hat I_{\T}$.

Given $\mu\in \M_{erg}$, we say that $\mu$ \emph{lifts to $\hat I$} if there exists an ergodic $\hat f$-invariant probability measure $\hat\mu$ on $\hat I$ such that $\hat\mu\circ\pi^{-1}=\mu$.  For $f\in \F$, if $\mu\in \M_{erg}$ and $\lambda(\mu)>0$ then $\mu$ lifts to $\hat I$, see \cite{Kellift, BrKel}.

For convenience later, we let $\iota:=\pi|_{D_0}^{-1}$.  Note that there is a natural distance function $d_{\hat I}$ within domains $D$ (but not between them) induced from the Euclidean metric on $I$.

\subsection{Inducing schemes}

We say that $(X,F,\tau)$ is an \emph{inducing scheme} for $(I,f)$ if
\begin{list}{$\bullet$}{\itemsep 0.2mm \topsep 0.2mm \itemindent -0mm \leftmargin=5mm}
\item $X$ is an interval containing a finite or countable
collection of disjoint intervals $X_i$ \st $F$ maps each $X_i$
diffeomorphically onto $X$, with bounded distortion (i.e. there
exists $K>0$ so that for all $i$ and $x,y\in X_i$, $1/K\le DF(x)/DF(y) \le K$);
\item $\tau|_{X_i} = \tau_i$ for some $\tau_i \in \N$ and $F|_{X_i} = f^{\tau_i}$.  If $x \notin \cup_iX_i$ then $\tau(x)=\infty$.
\end{list}
The function $\tau:\cup_i X_i \to \N$ is called the {\em inducing time}. It may happen that $\tau(x)$ is the first return time of $x$ to $X$, but
that is certainly not the general case.  For ease of notation, we will often write $(X,F)=(X,F,\tau)$.  In this paper we can always assume that every inducing scheme is uniformly expanding.

Given an inducing scheme $(X,F, \tau)$, we say that a measure $\mu_F$ is a \emph{lift} of $\mu$ if for all $\mu$-measurable subsets $A\subset I$,
\begin{equation} \mu(A) = \frac1{\int_X \tau \ d\mu_F} \sum_i \sum_{k = 0}^{\tau_i-1} \mu_F( X_i \cap f^{-k}(A)). \label{eq:lift}
\end{equation}
Conversely, given a measure $\mu_F$ for $(X,F)$, we say that
$\mu_F$ \emph{projects} to $\mu$ if \eqref{eq:lift} holds.
We denote $$(X,F)^\infty:=\left\{x\in X:\tau(F^k(x)) \text{ is defined for all }k\ge 0\right\}.$$
We call a measure
$\mu$  \emph{compatible to} the inducing scheme $(X,F,\tau)$ if
\begin{list}{$\bullet$}{\itemsep 0.2mm \topsep 0.2mm \itemindent -0mm \leftmargin=5mm}
\item $\mu(X)> 0$ and $\mu(X \setminus (X,F)^{\infty}) = 0$; and
\item there exists a measure $\mu_F$ which projects to $\mu$ by
\eqref{eq:lift}, and in particular $\int_X \tau \ d\mu_F <
\infty$.
\end{list}

For a potential $\phi:I \to \R$, we define the \emph{induced potential} $\Phi:X\to \R$ for an inducing scheme $(X,F,\tau)$ as $$\Phi(x):=S_{\tau(x)}\phi(x)=\phi(x)+\ldots+\phi\circ f^{\tau(x)-1}(x)$$ whenever $\tau(x)<\infty$.  We denote $\Phi_i:=\sup_{x\in X_i}\Phi(x)$.  Note that sometimes we will abuse notation and write $(X,F,\Phi)$ when we are particularly interested in the induced potential for the inducing scheme.
The following is known as Abramov's formula, see for example \cite{Zwei, PeSe}.

\begin{lemma}
Let $\mu_F$ be an ergodic invariant measure on  $(X,F,\tau)$ such that $\int\tau~d\mu_F<\infty$ and with projected measure $\mu$.
Then $h_{\mu_F}(F)=\left(\int\tau~d\mu_F\right)h_\mu(f)$.  Moreover, if $\phi:I\to \R$ is a potential, and $\Phi$ the corresponding induced potential, then $\int\Phi~d\mu_F = \left(\int\tau~d\mu_F\right)\int\phi~d\mu$.
\label{lem:abra}
\end{lemma}

Fixing $f$, we let
$$\M_+:=\{\mu\in \M_{erg}:\lambda_f(\mu)>0\}, \text{ and for } \eps>0, \ \M_\eps:=\{\mu\in \M_{erg}: h_\mu\ge\eps\}.$$

For a proof of the following result, see \cite[Theorem 3]{BTeqnat}.

\begin{theorem}
If $f\in \F$ and $\mu \in \M_+$, then there is an inducing scheme $(X,F,\tau)$
and a measure $\mu_F$ on $X$ \st $\int_X \tau \ d\mu_F < \infty$.
Here $\mu_F$ is the lifted measure of $\mu$ (i.e., $\mu$ and
$\mu_F$ are related by \eqref{eq:lift}).  Moreover,
$\overline{(X,F)^\infty}=X\cap\Omega$.

Conversely, if $(X,F,\tau)$ is an inducing scheme and $\mu_F$ an
ergodic $F$-invariant measure \st $\int_X \tau d\mu_F < \infty$, then $\mu_F$ projects to a measure $\mu \in \M_+$.
\label{thm:schemes}
\end{theorem}

The proof of the above theorem uses the theory of \cite[Section 3]{BrCMP}.  The main idea is that the Hofbauer extension can be used to produce inducing schemes.
We pick $\hat X\subset \hat I_{\T}$ and use a first return map to $\hat X$ to give the inducing scheme on $X:=\pi(\hat X)$.  We will always choose $X$ to be a cylinder in $\P_n$, for various values  of $n\in \N$.  As in \cite{BTeqnat}, sets $\hat X$, and thus the inducing schemes they give rise to, will be of two types.

\textbf{Type A:} The set $\hat X$ is an interval in a single domain $D\in \D_{\T}$.  Then for $x\in X$ there exists a unique $\hat x\in \hat X$ so that $\pi(\hat x)=x$.  Then $\tau(x)$ is defined as the first return time of $\hat x$ to $\hat X$.  We choose $\hat X$ so that $X\in \P_n$ for some $n$, and $\hat X$ is compactly contained in $D$.  These properties mean that $(X,F,\tau)$ is an inducing scheme which is extendible.  That is to say, letting $X'=\pi(D)$, for any domain $X_i$ of $(X,F)$ there is an extension of $f^{\tau_i}$ to $X_i'\supset X_i$ so that $f^{\tau_i}:X_i' \to X'$ is a homeomorphism.  By the distortion \cite[Theorem C(2)]{SVarg}, this means that $(X,F)$ has uniformly bounded distortion, with distortion constant depending on $\delta:=d_{\hat I}(\hat X, \bd D)$.

\textbf{Type B:}
We fix $\delta>0$ and some interval $X\in \P_n$ for some $n$.
We say that the interval $X'$ is a \emph{$\delta$-scaled neighbourhood of $X$} if, denoting the left and right components of $X' \sm X$ by $L$ and $R$ respectively, we have $|L|,|R|=\delta|X|$.
We fix such an $X'$ and let $\hat X = \sqcup \{ D \cap \pi^{-1}(X) : D \in \D_{\T}, \pi(D) \supset X'\}$.  Let $r_{\hat X}$ denote the first return time to $\hat X$. Given $x\in X$, for any $\hat x\in \hat X$ with $\pi(\hat x)=x$, we set $\tau(x)=r_{\hat X}(\hat x)$.   In \cite{BrCMP} it is shown that by the setup, this time is independent of the choice of $\hat x$ in $\pi|_{\hat X}^{-1}(x)$.  Also for each $X_i$ there exists $X_i'\supset X_i$ so that $f^{\tau_i}:X_i' \to X'$ is a homeomorphism, and so, again by the Koebe Lemma, $F$ has uniformly bounded distortion, with distortion constant depending on $\delta$.

We will need to deal with both kinds of inducing scheme since we want information on the tail behaviour, \ie the measure of $\{\tau\ge n\}$ for different measures.  As in Propositions~\ref{prop:eqgen} and \ref{prop:eqnat} below, for measures close to $\mu_\phi$ we have good tail behaviour for schemes of type A; and for measures close to the acip $\mu_{-\log|Df|}$ we have good tail behaviour for schemes of type B.  We would like to point out that any type A inducing time $\tau_1$ can be expressed as a power of a type B inducing time $\tau_2$, \ie $\tau_1= \tau_2^p$ where $p:X \to \N$.
Moreover, $\int p~d\mu_1<\infty$ for the induced measure $\mu_1$ for the type A inducing scheme.  This type of relation is considered by Zweim\"uller \cite{Zwei}.

\subsection{Method of proof}

The main difficulty in the proof of Theorem~\ref{thm:main spectrum} is to get an upper bound on the dimension spectrum in terms of $T_\phi$.  To do this, we show that there are inducing schemes which have sufficient multifractal information to give an upper bound on $\widetilde{\DS}_\phi$.  Then we can use Iommi's main theorem in \cite{Iom}, which gives upper bounds in terms of the $T$ for the inducing scheme.  It is the use of these inducing schemes which is the key to this paper.

We first show in Section~\ref{sec:param} that for a given range of $\alpha$ there are inducing schemes which are compatible to any measure $\mu$ which has $h_\mu+\int\psi_q~d\mu$ sufficiently large, where $q$ depends on $\alpha$.  In doing this we will give most of the theory of thermodynamic formalism needed in this paper.  For example, we show the existence of equilibrium states on $\K_\alpha$ which will turn out to have full dimension (these also give the lower bound for $\widetilde{\DS}_\phi$).

In Section~\ref{sec:lyap} we prove that for a set $A$, there is an inducing scheme that `sees' all points $x\in A$ with $\overline\lambda_f(x)$ bounded below, up to set of small Hausdorff dimension.  This means that we can fix inducing schemes which contain all the relevant measures, as above, and also contain the multifractal data.  Then in Section~\ref{sec:main thm} we prove Theorem~\ref{thm:main spectrum} and Proposition~\ref{prop:nonreg}.  In Section~\ref{sec:LE spec} we show how our results immediately give us information on the Lyapunov spectrum.  In the appendix we show that pointwise dimensions for induced measures and the original ones are the same, also extending our results to potentials in the class $SVI$.

\section{The range of parameters}
\label{sec:param}

In this section we determine what $U$ is in Theorem~\ref{thm:main spectrum}.  In order to do so, we must introduce most of the theory of the thermodynamical properties for inducing schemes required in this paper.  Firstly we show that if $\alpha(q)\in U$, then the equilibrium states for $\psi_q$ are forced to have positive entropy.  By Theorem~\ref{thm:schemes}, this ensures that the equilibrium states must be compatible to some inducing scheme, and thus we will be able to use Iommi's theory.

We let $$G_\eps(\phi):=\left\{q:\exists\delta<0 \text{ such that } \int\psi_q~d\mu>\delta \Rightarrow h_\mu>\eps\right\}.$$  The next lemma shows that most of the relevant parameters $q$ which we are interested in must lie in $G_\eps(\phi)$.

\begin{lemma}  Let $\phi:I \to \R$ be a potential satisfying \eqref{eq:range} and with $P(\phi)=0$.  Suppose that \eqref{eq:grow} holds for $f$.
There exist $\eps>0$, $q_1<1<q_2$ so that $(q_1,q_2)\subset G_\eps$.  If we take $\eps>0$ arbitrarily close to 0 then we can take $q_1$ arbitrarily close to 0.  If \eqref{eq:CE} holds then $[0,1]\subset (q_1, q_2)$.
\label{lem:entropy}
\end{lemma}

\begin{proof}
We first prove $q_1<1$.  By \eqref{eq:range} and $P(\phi)=0$, we must have $\phi<0$:
$$
0=P(\phi) \ge \htop(f)+\int \phi~d\mu_{-\htop(f)}\ge \htop(f) + \inf\phi > \sup\phi
$$ where $\mu_{-\htop(f)}$ denotes the measure of maximal entropy (for more details of this measure, see Section~\ref{sec:LE spec}).
Let $q_1$ be any value in $(0,1]$. Then suppose that for some $\delta<0$, a measure $\mu\in \M_{erg}$ has $$h_\mu+\int-T_\phi(q_1)\log|Df|+q_1\phi~d\mu>\delta$$ for $T_\phi$ as in \eqref{eq:Tphi}.  Recall that by \cite{Prz}, $\lambda(\mu)\ge 0$ since we excluded the possibility of attracting cycles for maps $f\in \F$.   Then
$$h_\mu>\delta+\int T_\phi(q_1)\log|Df|-q_1\phi~d\mu \ge \delta+q_1|\sup\phi|>0.$$
If $q_1$ was chosen very close to 0 then $\eps>0$ must be chosen small too.

We can similarly show that $q_2>1$, the only difference in this case being that $q>1$ implies that $T_\phi(q)<0$.  So we can take $\mu$ as above and obtain
$$h_\mu>\delta+\int T_\phi(q)\log|Df|-q\phi~d\mu \ge \delta+\int T_\phi(q)\log|Df|~d\mu+q|\sup\phi|.$$
Since $T_\phi(q)$ is close to 0 for $q$ close to 1 and since $\int \log|Df|~d\mu<\log\sup_{x\in I} |Df(x)|$, for $q_2>1$ close to 1, the above must be strictly positive.

Suppose now that \eqref{eq:CE} holds.  Then by \cite{BrvS}, there exists $\eta>0$ so that any invariant measure $\mu\in \M_{erg}$ must have $\lambda_f(\mu)>\eta$.  So if
$h_\mu+\int-T_\phi(q)\log|Df|+q\phi~d\mu>\delta$, then
$$h_\mu>\delta+\int T_\phi(q)\log|Df|-q\phi~d\mu \ge \delta+T_\phi(q)\eta+q|\sup\phi|.$$  For $q$ close to 0, $T_\phi(q)$ must be close to 1, so we can choose $\delta$ and $q_1<0$ so that the lemma holds.
\end{proof}

\textbf{The sets $\mathbf{Cover(\eps)}$ and $\mathbf{SCover(\eps)}$:}
Let $\eps>0$.  By \cite[Remark 6]{BTeqnat} there exists $\eta>0$ and a compact set $\hat E\subset \hat I_{\T}$ so that $\mu\in \M_\eps$ implies that $\hat\mu(\hat E)>\eta$.  Moreover $\hat E$ can be taken inside $\hat I_R\sm B_\delta(\bd\hat I)$ for some $R\in \N$ and $\delta>0$. (Here $B_\delta(\bd\hat I)$ is a $\delta$-neighbourhood of $\bd\hat I$ with respect to the distance function $d_{\hat I}$).  As in \cite[Section 4.2]{BTeqnat}, $\hat E$ can be covered with sets $\hat X_1, \ldots, \hat X_n$ so that each $\hat X_k$ acts as the set which gives the inducing schemes $(X_k,F_k)$ (where $X_k=\pi(\hat X_k)$) as in Theorem~\ref{thm:schemes}.  We will suppose that these sets are either all of type A, or all of type B.  This means that any $\mu\in \M_\eps$ must be compatible to at least one of $(X_k, F_k)$.  We denote $Cover^A(\eps)=\{\hat X_1, \ldots, \hat X_n\}$ and the corresponding set of schemes by $SCover^A(\eps)$ if we are dealing with type A inducing schemes.  Similarly we use $Cover^B(\eps)$ and  $SCover^B(\eps)$ for type B inducing schemes.  If a result applies to both schemes of types then we omit the superscript.

We let $\{X_{k,i}\}_{i}$ denote the domains of the inducing scheme $(X_k,F_k)$ and we denote the value of $\tau_k$ on $X_{k,i}$ by $\tau_{k,i}$.  Given $(X_k,F_k, \tau_k)$, we let $\Psi_{q,k}$ denote the induced potential for $\psi_q$.%, and let $\Psi_{q,k,i}:=\sup_{x\in X_{k,i}}\Psi_{q,k}(x)$.

From this setup, given $q\in G_\eps(\phi)$ there must exist a sequence of measures $\{\mu_n\}_n\subset \M_\eps$ and a scheme $(X_k, F_k)$ so that $h_{\mu_n}+\int\psi_q~d\mu_n \to P(\psi_q)=0$ and $\mu_n$ are all compatible to $(X_k,F_k)$.  Later this fact will allow us to use \cite[Proposition 1]{BTeqnat} to study equilibrium states for $\psi_q$.

If $\upsilon:I \to \R$ is some potential and $(X,F)$ is an inducing scheme with induced potential $\Upsilon:X\to \R$, we let $\Upsilon_i:=\sup_{x\in X_i}\Upsilon(x)$.  We define
the \emph{$k$th variation} as $$V_k(\Upsilon):=\sup_{\c_k\in \P_k}\{|\Upsilon(x)-\Upsilon(y)|:x,y\in \c_k\}.$$  We say that $\Upsilon$ is \emph{locally H\"older continuous} if there exists $\alpha>0$ so that  $V_k(\Upsilon)=O(e^{-\alpha n})$.  We let
\begin{equation} Z_0(\Upsilon):=\sum_ie^{\Upsilon_i}, \text{ and } Z_0^*(\Upsilon):=\sum_i\tau_ie^{\Upsilon_i}.\label{eq:Z0}\end{equation}  As in \cite{SaBIP}, if $\Upsilon$ is locally H\"older continuous, then $Z_0(\Upsilon)<\infty$ implies $P(\Upsilon)<\infty$.

We say that a measure $\mu$ satisfies the \emph{Gibbs property with constant $P\in \R$} for $(X,F,\Upsilon)$ if there exists $K_\Phi, P\in \R$ so that $$\frac1{K_\Phi}\le \frac{\mu(\c_n)}{e^{S_n\Upsilon(x)-nP}} \le K_\Phi$$
for every $n$-cylinder $\c_n$ and all $x\in \c_n$.

The following is the main result of \cite{BTeqgen} (in fact it is proved for a larger class of potentials there).

\begin{proposition}
Given $f\in \F$ satisfying \eqref{eq:grow} and $\phi:I\to \R$ a H\"older potential satisfying \eqref{eq:range} and with $P(\phi)=0$, then for any $\eps>0$ and any $(X,F)\in SCover^A(\eps)$ with induced potential $\Phi$:
\begin{enumerate}[(a)]
\item There exists $\beta_\Phi>0$ such that $\sum_{\tau_i=n}e^{\Phi_i} =O(e^{-n\beta_\Phi})$;
\item $\Phi$ is locally H\"older continuous and $P(\Phi)=0$;
\item There exists a unique $\Phi$-conformal measure $m_\Phi$, and a unique equilibrium state $\mu_\Phi$ for $(X,F,\Phi)$.
\item There exists $C_\Phi$ so that $\frac1{C_\Phi}\le \frac{d\mu_\Phi}{dm_\Phi}\le C_\Phi$;
\item There exists an equilibrium state $\mu_\phi$ for $(I,f,\phi)$;
\item The map $t\mapsto P(t\phi)$ is analytic for $t\in  \left(\frac{-\htop(f)}{\sup\phi-\inf\phi}, \frac{\htop(f)}{\sup\phi-\inf\phi}\right)$.
\end{enumerate}
\label{prop:eqgen}
\end{proposition}

The existence of the equilibrium state $\mu_\phi$ under even weaker conditions than these was proved by Keller \cite{kellholder}.  However, we need all of the properties above to complete our analysis of the dimension spectrum of $\mu_\phi$.

The following is proved in \cite{BTeqnat}.  For the same result for unimodal maps satisfying \eqref{eq:CE} see \cite{BrKel}, which used tools from \cite{KelNow}.

\begin{proposition}
Suppose that $f\in \F$ satisfies \eqref{eq:poly} and let $$\psi(x)=\psi_t(x):=-t\log|Df(x)|-P(-t\log|Df(x)|).$$  Then there exists $t_0<1$ such that for any $t\in (t_0,1)$ there is $\eps=\eps(t)>0$ so that for any $(X,F)\in SCover^B(\eps)$ with induced potential $\Psi$:
\begin{enumerate}[(a)]
\item There exists $\beta_{DF}>0$ such that $\sum_{\tau_i=n}e^{\Psi_i} =O(e^{-n\beta_{DF}})$;
\item $\Psi$ is locally H\"older continuous and $P(\Psi)=0$;
\item There exists a unique $\Psi$-conformal measure $m_\Psi$, and a unique equilibrium state $\mu_\Psi$  for $(X,F,\Psi)$;
\item There exists $C_\Psi$ so that $\frac1{C_\Psi}\le \frac{d\mu_\Psi}{dm_\Psi}\le C_\Psi$;
\item There exists an equilibrium state $\mu_\psi$ for $(I,f,\psi)$ and thus for $(I, f, -\log|Df|)$;
\item The map $t\mapsto P(-t\log|Df|)$ is analytic in $(t_0,1)$.
\end{enumerate}
\label{prop:eqnat}
\end{proposition}

If $f\in \F$ satisfies \eqref{eq:CE}, then this proposition can be extended so that $t$ can be taken in a two-sided neighbourhood of 1.

In Proposition~\ref{prop:eqgen} both $m_\Phi$ and $\mu_\Phi$ satisfy the Gibbs property, and in Proposition~\ref{prop:eqnat} both $m_\Psi$ and $\mu_\Psi$ satisfy the Gibbs property: in all these cases, the Gibbs constant $P$ is 0.  By the Gibbs property, part (a) of Proposition~\ref{prop:eqgen} and \ref{prop:eqnat} imply that $\mu_\Phi(\{\tau=n\})$ and $\mu_\Psi(\{\tau=n\})$ respectively decay exponentially.  These systems are referred to as having \emph{exponential tails}.

One consequence of the first item in both of these propositions, as noted in \cite[Theorem 10]{BTeqgen} and \cite[Theorem 5]{BTeqnat}, is that we can consider combinations of the potentials above: $x\mapsto -t\log|Df(x)|+s\phi(x) -P(-t\log|Df|+s\phi)$.  We can derive the same results for this potential for $t$ close to 1 and $s$ sufficiently close to 0, or alternatively for $s$ close to 1 and $t$ sufficiently close to 0.  Note that by \cite{KelNow, BrKel} this can also be shown in the setting of unimodal maps satisfying \eqref{eq:CE} with potentials $\phi$ of bounded variation.

If $(X,F)$ is an inducing scheme with induced potential $\Phi:X\to \R$, we define
\iffalse
$$PB_\eps(\Phi):=\left\{q\in G_\eps(\phi): \exists \delta>0 \text{ s.t. } \sum_i\tau_ie^{\Psi_{q,i}+\delta\tau_{i}}<\infty\right\}.$$
\fi
$$PB_\eps(\Phi):=\left\{q\in G_\eps(\phi): \exists \delta>0 \text{ s.t. } Z_0^*(\Psi_{q}+\tau\delta)<\infty\right\}.$$

\begin{lemma}   For $(X_k,F_k)\in SCover(\eps)$, if $q\in PB_\eps(\Phi_k)$ then $P(\Psi_{q,k})=0$.  Moreover, there is an equilibrium state $\mu_{\Psi_{q,k}}$ for $(X_k,F_k,\Psi_{q,k})$  and the corresponding projected equilibrium state $\mu_{\psi_q}$ is compatible to any $(X_j,F_j)\in SCover(\eps)$.
\label{lem:PB Scover}
\end{lemma}

In this lemma, $SCover(\eps)$ can be $SCover^A(\eps)$ or $SCover^B(\eps)$.  Note that by \cite[Proposition 1]{BTeqnat}, if for any $(X,F)\in SCover(\eps)$ and $q\in PB_\eps(\Phi)$, then there exists an equilibrium state $\mu_{\Psi_q}$ for $(X,F,\Psi_q)$, as well as a unique equilibrium state $\mu_{\psi_q}$ for $(I,f,\psi_q)$.

\begin{proof}
Firstly we have $P(\Psi_{q,k})=0$ for the inducing scheme $(X_k, F_k)$ by Case 3 of \cite[Proposition 1]{BTeqnat}.  Secondly we can replace $(X_k,F_k)$ with any inducing scheme $(X_j,F_j)\in SCover(\eps)$ by \cite[Lemma 9]{BTeqnat}.
\end{proof}

This lemma means that if $q\in PB_\eps(\Phi_k)$ for $(X_k,F_k)\in SCover^A(\eps)$, then $q\in PB_\eps(\Phi_j)$ for any $(X_j,F_j)\in SCover^A(\eps)$.  Therefore, we can denote this set of $q$ by  $PB_\eps^A(\phi)$.  Since the same argument holds for inducing schemes of type B, we can analogously define the set $PB_\eps^B(\phi)$.  Note that $\eps'<\eps$ implies $PB_{\eps'}(\phi)\supset PB_\eps(\phi)$.  We define $$PB(\phi):=\cup_{\eps>0}PB_\eps(\phi).$$

\begin{remark}
The structure of inducing schemes here means that we could just fix a single inducing scheme which has all the required thermodynamic properties in this section.  However, in Section~\ref{sec:lyap} we  need to consider all the inducing schemes here in order to investigate the dimension spectrum.
\end{remark}

In \cite{Iom}, the following conditions are given.

$$q^*:=\inf\{q:\text{there exists } t\in \R\text{ such that } P(-t\log|DF|+q\Phi)\le 0\}.$$

\begin{equation*}
T_\Phi(q):=\begin{cases}
\inf\{t\in \R: P(-t\log|DF|+q\Phi)\le 0\} & \text{if } q\ge q^*,\\
\infty & \text{if } q<q^*.
\end{cases}
\end{equation*}

The following is the main result of \cite[Theorem 4.1]{Iom}.  We can apply it to our schemes $(X,F)$ since they can be seen as the full shift on countably many symbols $(\Sigma, \sigma)$.  In applying this theorem, we choose the metric $d_\Sigma$ on $\Sigma$ to be compatible with the Euclidean metric on $X$.

\begin{theorem}
Suppose that $(\Sigma,\sigma)$ is the full shift on countably many symbols and $\Phi:\Sigma\to \R$ is locally H\"older continuous. The dimension spectrum $\alpha\mapsto \DS_\Phi(\alpha)$ is the Legendre transform of $q\mapsto T_\Phi(q)$. \label{thm:Iommi}
\end{theorem}

If we know that an inducing scheme has sufficiently high, but not infinite, pressure for the potential $\Psi_q$ then the measures we are interested in are all compatible to this inducing scheme.  This leads to $T_\Phi$ defined above being equal to $T_\phi$ as defined in \eqref{eq:Tphi}, as in the following proposition.

\begin{proposition}
Suppose that $f\in \F$ is a map satisfying \eqref{eq:grow} and $\phi:I \to \R$ is a H\"older potential satisfying \eqref{eq:range}.  Let $\eps>0$.  For all $q\in PB_\eps^A(\phi)$, if $(X,F)\in SCover^A(\eps)$ with induced potential $\Phi$, then $T_\Phi(q)=T_\phi(q)$.
Similarly for type B inducing schemes.

Moreover,
\begin{itemize}
\item[(a)] there exists $\eps>0$ and $q_0<1<q_1$ so that $(q_0, q_1)\subset PB_\eps^A(\phi)$;

\item[(b)] if $f$ satisfies \eqref{eq:poly}, then for all $\eps>0$ there exist $0<q_2<q_3$ so that $(q_2, q_3)\subset PB_\eps^B(\phi)$ (taking $\eps$ small, $q_2$ can be taken arbitrarily close to 0);
\item[(c)] if $f$ satisfies \eqref{eq:CE}, for all $\eps>0$ there exist $q_2<0<q_3$ so that $(q_2, q_3)\subset PB_\eps^B(\phi)$.
\end{itemize}
\label{prop:q nhd}
\end{proposition}

In this proof, and later in the paper, given a set $A$ and a function $g:A\to \R$ we let $$|g|_\infty=\sup_{x\in A}|g(x)|.$$

\begin{proof}
By Lemma~\ref{lem:PB Scover}, for $q\in PB_\eps(\phi)$, and any $(X,F)\in SCover(\eps)$, $P(\Psi_q)=0$.  The Abramov formula in Lemma~\ref{lem:abra} implies that
\begin{align*}
0&=h_{\mu_{\psi_q}}(f)+\int-T_\phi(q)\log|Df|+q\phi~d\mu_{\psi_q}\\
&= \left(\frac1{\int\tau~d\mu_{\Psi_q}}\right)
\left(h_{\mu_{\Psi_q}}(F)+\int-T_\phi(q)\log|DF|+q\Phi~d\mu_{\Psi_q}\right) \end{align*} and hence $T_\Phi(q)\le T_\phi(q)$ on $PB_\eps(\phi)$.  Since $\log|DF|$ is uniformly positive, we also know that $t\mapsto P(-t\log|DF|+q\Phi)$ is strictly decreasing in $t$ and hence $T_\Phi(q)=T_\phi(q)$ on $PB_\eps(\phi)$.

By Lemma~\ref{lem:PB Scover}, for $\eps>0$, in order to check if $q\in PB_\eps(\phi)$ and thus prove (a), (b) and (c), we only need to check if $q\in PB_\eps(\Phi)$ for one scheme $(X,F)\in SCover(\eps)$. We will show that the estimate for $Z_0^*(\Psi_q)$ is a sum of exponentially decaying terms, which is enough to show that there exists $\delta>0$ so that $Z_0^*(\Psi_q+\delta\tau)<\infty$.

As in the proof of Lemma~\ref{lem:entropy}, \eqref{eq:range} and $P(\phi)=0$ imply that $\phi<0$.
Recall that by definition $P(-T_\phi(q)\log|Df|+q\phi)=0$.  Given $(X,F)\in SCover(\eps)$, by the local H\"older continuity of every $\Psi_q$, there exists $C>0$ such that for $Z_0^*$ as in \eqref{eq:Z0}, $$Z_0^*(\Psi_q):=\sum_i \tau_ie^{-T_\phi(q)\log|DF_i|+q\Phi_i}\le C
\sum_nn\sum_{\tau_i=n} |X_i|^{T_\phi(q)}e^{q\Phi_i}.$$
We will first assume only that $f$ satisfies \eqref{eq:grow} and that $q$ is close to 1.  In this case we work with inducing schemes of type A.  By Proposition~\ref{prop:eqgen}, there exists $\beta_\Phi>0$ so that  $\sum_{\tau_i=n}e^{\Phi_i}=O(e^{-n\beta_\Phi})$.

\textbf{Case 1:} $q$ near 1 and $q>1$.  In this case $T_\phi(q)<0$.  Since $|X_i|\ge (|Df|_\infty)^{-\tau_i}$,
$$Z_0^*(\Psi_q)\le C
\sum_nn(|Df|_\infty)^{-n|T_\phi(q)|} \sum_{\tau_i=n}e^{q\Phi_i}\le
C' \sum_nn(|Df|_\infty)^{-n|T_\phi(q)|} e^{-nq\beta_\Phi}.$$
Since for $q$ near to 1, $T_\phi(q)$ is close to 0, the terms on the right decay exponentially, proving the existence of $q_1>1$ in part (a).

\textbf{Case 2:} $q$ near 1 and $q<1$.  In this case $T_\phi(q)>0$.  By the H\"older inequality there exists $C'>0$ such that
\begin{align*}
Z_0^*(\Psi_q)&\le C\sum_nn\left(\sum_{\tau_i=n} e^{\Phi_i}\right)^q \left(\sum_{\tau_i=n} |X_i|^{\frac{T_\phi(q)}{1-q}}\right)^{1-q}\\
&\le C'\sum_n n e^{-qn\beta_\Phi} \left(\sum_{\tau_i=n}|X_i|^{\frac{T_\phi(q)}{1-q}}\right)^{1-q}.\end{align*}

\textbf{Case 2(a):} $\frac{T_\phi(q)}{1-q} \ge 1$.
In this case obviously $Z_0^*(\Psi_q)$ can be estimated by exponentially decaying terms.  (In fact, it is not too hard to show that this case is empty, but there is no need to give the details here.)

\textbf{Case 2(b):} $\frac{T_\phi(q)}{1-q}<1$.  Here the term we need to control is, by the H\"older inequality
$$\left(\sum_{\tau_i=n}|X_i|^{\frac{T_\phi(q)}{1-q}}\right)^{1-q}\le \left[\left(\sum_{\tau_i=n}|X_i|\right)^{\frac{T_\phi(q)}{1-q}}
\#\{\tau_i=n\}^{1-\left(\frac{T_\phi(q)}{1-q}\right)}\right]^{1-q}.$$
We have $$\left(\#\{\tau_i=n\}^{1-\left(\frac{T_\phi(q)}{1-q}\right)}\right)^{1-q}= \#\{\tau_i=n\}^{1-q -T_\phi(q)}.$$  As explained in \cite{BTeqnat}, for any $\eta>0$ there exists $C_\eta>0$ such that $\#\{\tau_i=n\}\le C_\eta e^{n(\htop(f)+\eta)}$.  Since we also know that for $q$ close to 1, $1-q -T_\phi(q)$ is close to 0, the terms $e^{-nq\beta_\Phi}$ dominate the estimate for $Z_0^*(\Psi_q)$, which completes the proof of part (a) of the proposition.

Next we assume that $f$ satisfies \eqref{eq:poly} and $q>0$ is close to 0.  In this case we work with inducing schemes of type B.

\textbf{Case 3:} $q$ near 0 and $q>0$.  In this case $T_\phi(q)<1$.  By \cite[Proposition 3]{BTeqnat}, if $t$ is close to 1 then $\sum_{\tau_i=n} |X_i|^{t}$ is uniformly bounded.  Thus, as in Case 2,
$$Z_0^*(\Psi_q)\le C\sum_nn\left(\sum_{\tau_i=n} e^{\Phi_i}\right)^q \left(\sum_{\tau_i=n} |X_i|^{\frac{T_\phi(q)}{1-q}}\right)^{1-q}=O\left(n\sum_{\tau_i=n} e^{\Phi_i}\right)^q.$$
As in Case 2, there exists $\beta_\Phi>0$ so that $\mu_{\Phi}\{\tau=n\}=O(e^{-n\beta_\Phi})$, which implies $Z_0^*(\Psi_q)$ can be estimated by exponentially decaying terms, proving (b).

\textbf{Case 4:} $q$ near 0 and $q<0$.  This can only be considered when $f$ satisfies \eqref{eq:CE}.  In this case $T_\phi(q)>1$.   Note that by Proposition~\ref{prop:eqnat} there exists $\beta_{DF}>0$ so that $\mu_{-\log|DF|}\{\tau=n\}=O(e^{-n\beta_{DF}})$.  Thus,
$$Z_0^*(\Psi_q)\le C\sum_nne^{qn\inf\phi}\left(\sum_{\tau_i=n} |X_i|\right)^{T_\phi(q)} =O\left(\sum_nne^{n(q\inf\phi-T_\phi(q)\beta_{DF})}\right).$$
For $q$ close to 0 we have $q\inf\phi-T_\phi(q)\beta_{DF}<0$ and so
$Z_0^*(\Psi_q)$ can be estimated by exponentially decaying terms, proving (c).
\iffalse
We claim that the above computations also give us that $\int\tau~d\mu_{\Psi_q}<\infty$ for the domains of $q$ considered.  This is because $\int\tau~d\mu_{\Psi_q}=\sum_i \tau_i \mu_{\Psi_q}(X_i)$, which by the Gibbs property of $\Psi_q$ is bounded above by $O\left(\sum_i \tau_i e^{\sup_{x\in X_i}\Psi_q(x)}\right)$.  Since the above computations all give exponential decay of the terms $\sum_{\tau_i=n}e^{\sup_{x\in X_i}\Psi_q(x)}$, this quantity is bounded, as required.\fi
\end{proof}

\section{Inducing schemes see most points with positive Lyapunov Exponent}
\label{sec:lyap}

The purpose of this section is to show that if we are only interested in those sets for which the Lyapunov exponent is bounded away from 0, then there are inducing schemes which contain all the multifractal data for these sets.  This is the content of the following proposition.

\begin{proposition}
For all $\lambda,s>0$ there exist $\eps=\eps(\lambda,s)>0$, a set $\oG_\lambda' \subset  \oG_\lambda$, and an inducing scheme $(X,F)\in SCover^A(\eps)$ so that $\dim_H(\oG_\lambda\sm\oG_\lambda')\le s$ and for all $x\in \oG_\lambda'$ there exists $k\ge 0$ so that $f^k(x)\in (X,F)^\infty$.  There is also an inducing scheme in $SCover^B(\eps)$ with the same property.
\label{prop:LE}
\end{proposition}

By the structure of the inducing schemes outlined above, we can replace $\eps$ with any $\eps'\in (0,\eps)$.
This means that if there is a set $A\subset I$ and $\lambda>0$ so that $\dim_H(A\cap\oG_\lambda)>0$ then there is an inducing scheme $(X,F)$ so that $\dim_H(A\cap\oG_\lambda\cap (X,F)^\infty)= \dim_H(A\cap\oG_\lambda)$.  Hence the multifractal information for $A\cap\oG_\lambda$ can be found using $(X,F)$.
We remark that by Lemma~\ref{lem:PB Scover}, for $\lambda>0$ and $q\in PB(\phi)$, if $\dim_H(\K_\alpha\cap\oG_\lambda)>0$ then we can fix an inducing scheme $(X,F)$ such that $$\dim_H(\K_\alpha\cap\oG_\lambda\cap(X,F)^\infty) =
\dim_H(\K_\alpha\cap\oG_\lambda).$$

For the proof of Proposition~\ref{prop:LE} we will need two lemmas.

Partly for completeness and partly in order to fix notation, we recall the definition of Hausdorff measure and dimension.  For $E\subset I$ and $s,\delta>0$, we let $$H_\delta^s(E):=\inf\left\{\sum_i\diam(A_i)^s\right\}$$
where the infimum is taken over collections $\{A_i\}_i$ which cover $E$ and with $\diam(A_i)<\delta$.  Then the $s$-Hausdorff measure of $E$ is defined as $H^s(E):=\limsup_{\delta\to 0}H_\delta^s(E)$.    The Hausdorff dimension is then $\dim_H(E):=\sup\{s:H^s(E)=\infty\}$.

\begin{lemma}
For all $\lambda,s>0$ there exists $\eta>0$, $R\in \N$ and $\oG_\lambda'\subset \oG_\lambda$ so that  $\dim_H(\oG_\lambda \sm \oG_\lambda')\le s$, and $x\in \oG_\lambda'$  implies
$$\limsup_k\frac1k\#\{1\le k\le n:\hat f^k(\iota(x))\in \hat I_R\}>\eta.$$
\label{lem:big scale}
\end{lemma}

\iffalse
The main result of \cite{BrvS} shows that in the case that $f$ satisfies \eqref{eq:CE}, the Collet Eckmann condition, then there exists $\lambda>0$ so that $\overline\lambda_f(x)\ge\lambda$ for all $x\notin \cup_{n\ge 0}f^{-n}(\crit)$.  This exceptional set is countable, so has Hausdorff dimension 0.  This means that the lemma holds on all but a zero .....
\fi

\textbf{Note on the proof:}
It is important here that we can prove this lemma for $\oG_\lambda$ rather than $\uG_\lambda$.  Otherwise Proposition~\ref{prop:LE} and, for example, our main corollaries would not hold.  We would like to briefly discuss why we can prove this result for $\oG_\lambda$ rather than $\uG_\lambda$.  The argument we use in the proof is similar to arguments which show that under some condition on pointwise Lyapunov exponents for $m$-almost every point, there is an invariant measure absolutely continuous with respect to $m$.  Here $m$ is usually a conformal measure.  For example in \cite[Theorem 4]{BTholo} we showed that if $m(\uG_\lambda)>0$ for a conformal measure $m$ then `most points' spend a positive frequency of their orbit in a compact part of the Hofbauer extension, and hence there is an absolutely continuous invariant measure $\mu\ll m$.  In that case it was convenient to use $\uG_\lambda$ rather than $\oG_\lambda$.  In \cite{Kelthm}, and in a similar proof in \cite[Theorem V.3.2]{MSbook}, $m$ is Lebesgue measure and the ergodicity of $m$ is used to allow them to weaken assumptions and to consider $\oG_\lambda$ instead.  In our case here, we cannot use a property like ergodicity, but on the other hand we do not need points to enter a compact part of the extension with positive frequency (which is essentially  what is required in all the above cases), but simply infinitely often.  Hence we can use $\oG_\lambda$ instead.

For the proof of the lemma we will need the following result from \cite[Theorem 4]{BRSS}.  Here $m$ denotes Lebesgue measure, and as above $\ell_{max}(f)$ is the maximal critical order of all critical points of $f$.

\begin{proposition}
If $f\in \F$ satisfies \eqref{eq:grow} then there exists $K_0>0$ so that for any Borel set $A$,
$$m(f^{-n}(A))\le K_0m(A)^{\frac1{2\ell_{max}(f)}}.$$
\label{prop:BRSS}
\end{proposition}

\begin{remark}
For $f\in\F$, such a theorem holds whenever there is an acip $\mu$ with density $\rho=\frac{d\mu}{dm}\in L^{1+\delta}$ for $\delta>0$.  Standard arguments show that transitivity implies that there exists $\eps>0$ such that $\rho \ge\eps$.  Then \begin{align*}
m(f^{-n}(A)) &\le \frac1\eps \mu(f^{-n}(A))=
\frac1\eps \mu(A)= \frac1\eps \int 1_A\rho~dm\\
& \le \frac1\eps\left(\int 1_A~dm\right)^{\frac\delta{1+\delta}}\left(\int\rho^{1+\delta}~dm\right)^{\frac1{1+\delta}} \le Cm(A)^{\frac\delta{1+\delta}},
\end{align*}
for some $C>0$.
\end{remark}

\begin{proof}[Proof of Lemma~\ref{lem:big scale}]  For this proof we use ideas of \cite{Kellift}, see also \cite{BTholo}.  We also use the notation $|\cdot|$ to denote the length of a connected in interval.  We suppose that $\dim_H(\oG_\lambda)>0$, otherwise there is nothing to prove.  We fix $s\in (0,\dim_H(\oG_\lambda))$.  Throughout this proof, we write $\ell_{max}=\ell_{max}(f)$.

For $\gamma\ge 0$ and $n\in \N$, let $LG_\gamma^n:=\{x:|Df^n(x)|\ge e^{\gamma n}\}$.

For $x\in I$, we define $$\freq(R,\eta,n):= \left\{x:\frac1n\#\left\{0\le k<n:\hat f^k(\iota(x)) \in \hat I_R\right\}\le \eta\right\}$$ and
$$\freq(R,\eta):=\left\{y:\limsup_k\frac1k\#\left\{1\le k\le n:\hat f^k(\iota(y))\in \hat I_R\right\}<\eta\right\}.$$
For $\lambda_0\in(0, \lambda)$, $R, n \ge 1$ and $\eta>0$ we consider the set
$$E_{\lambda_0,R,n}(\eta):=LG_{\lambda_0}^n\cap\freq(R,\eta,n).$$

If $x\in \oG_\lambda\cap \freq(R,\eta)$ then there exists arbitrarily large $n\in \N$ so that $|Df^n(x)|\ge e^{\lambda_0 n}$, and $x\in \freq(R,\eta,n)$.  Hence
$$\freq(R,\eta)\cap\oG_\lambda \subset \bigcap_k\bigcup_{n\ge k} E_{\lambda_0,R,n}(\eta).$$  This means we can estimate the Hausdorff dimension of $\freq(R,\eta)\cap \oG_\lambda$ through estimates on $\dim_H(E_{\lambda_0,R,n}(\eta))$.

We let $\P\negmedspace_{E,n}$ denote the collection of cylinder sets of $\P_n$
which intersect $E_{\lambda_0,R,n}(\eta)$.
We will compute
$H_\delta^s(E_{\lambda_0,R,n}(\eta))$ using the natural structure of the dynamical cylinders $\P_n$.  First note that by \cite[Corollary 1]{Hpwise} (see also, for example, the proof of \cite[Theorem 4]{BTholo}), for all $\gamma>0$ there exist $R \ge 1$ and
$\eta>0$ so that $\#\P\negmedspace_{E,n}\le e^{\gamma n}$ for all large $n$.  In \cite{BTholo} this type of estimate was sufficient to show that conformal measure `lifted' to the Hofbauer extension.  The Hausdorff measure is more difficult to handle, since distortion causes more problems.  Here we use an argument of \cite{BTintstat} to deal with the distortion.  We will make some conditions on $\gamma$, depending on $s$ and $\lambda$ below.

Let $n(\delta)\in \N$ be so that $n\ge n(\delta)$ implies $|\c_n|<\delta$ for all $\c_n\in \P_n$.

We choose any $\gamma\in (0, \lambda s/16\ell_{max}^2)$ and
%$\gamma\in (0, \frac{\lambda s}{16\ell_{max}^2})$ and
$\theta:=4\gamma\ell_{max}^2/s$.
For $x\in I$, let $$V_n[x]:=\left\{y\in \c_n[x]:d(f^n(y),\bd f^n(\c_n[x]))<e^{-\theta n}|f^n(\c_n[x])|\right\}.$$
For a point $x\in E_{\lambda_0,R,n}$, we say that $x$ is in Case 1 if $x\in V_n[x]$, and in Case 2 otherwise.  We consider the measure of points in these different sets separately.

\textbf{Case 1:} For $x\in I$, we denote the part of $f^n(\c_n[x])$ which lies within $e^{-\theta n}|f^n(\c_n[x])|$ of the boundary of $f^n(\c_n[x])$ by $Bd_n[x]$.  We will estimate the Lebesgue measure of the pullback $f^{-n}(Bd_n[x])$.  Note that this set consists of more than just the pair of connected components $\c_n[x]\cap V_n[x]$.

Clearly, $m(Bd_n[x]) \le 2 e^{-\theta n}m(f^n(\c_n[x]))$.  Hence from  Proposition~\ref{prop:BRSS}, we have the (rather rough) estimate
\begin{align*}
m(V_n[x]) &\le m(f^{-n}(Bd_n[x]))
\le K_0\left[ 2 e^{-\theta n}m(f^n(\c_n[x]))\right]^{\frac1{2\ell_{max}^2}}\\
&\le 2K_0 e^{-\frac{\theta n}{2\ell_{max}^2}}=2K_0 e^{-\frac{2\gamma n}{s}}.
\end{align*}

\textbf{Case 2:} Let $\tilde \c_n[x]:=\c_n[x]\sm V_n[x]$.  As in the proof of \cite[Lemma 15]{BTintstat}, the intermediate value theorem and the Koebe lemma allow us to estimate $$\frac{|\tilde \c_n[x]|}{|f^n(\tilde \c_n[x])|}\le \left(\frac{1+e^{-n\theta}}{e^{-n\theta}}\right)^2\frac1{|Df^n(x)|}.$$
Hence for all large $n$, $$|\tilde\c_n[x]|\le 2e^{2\theta n}e^{-\lambda n}.$$
By our choice of $\gamma$,
$$|\tilde\c_n[x]|\le 2e^{-n\frac{\lambda}2}.$$

If we assume that $n\ge n(\delta)$, the sets $V_n[x]\subset \c_n[x]\in \P\negmedspace_{E,n}$ in Case 1 and $\tilde\c_n[x]\subset \c_n[x]\in \P\negmedspace_{E,n}$ in Case 2 form a $\delta$-cover of  $E_{\lambda_0,R,n}(\eta)$.  This implies that for $n$ large,
$$H_\delta^s(E_{\lambda_0,R,n}(\eta)) \le 4e^{\gamma n}(e^{-n\frac{\lambda s}2}+K_0 e^{-2\gamma n}).$$ By our choice of $\gamma$, this is uniformly bounded in $n$.  Since we can make the above estimate for all small $\delta$, we get that $$\dim_H\left( \oG_\lambda\cap\freq(R,\eta)\right)\le s.$$ So the set $\oG_\lambda':=\oG_\lambda \sm \freq(R,\eta)$ has the required property.
\end{proof}

Let $\{\eps_n\}_n$ be a positive sequence decreasing to 0 and let $B_n:=B_{\eps_n}(\bd \hat I)$, where we use the distance function $d_{\hat I}$ as described in Section~\ref{subsec:Hofbauer}.

\begin{lemma} For any $R\in \N$ and $\eta>0$, there exists $N(R,\eta)\in \N$ so that for $x\in I$, if $$\limsup_k\frac1k\#\left\{1\le j\le k:\hat f^j(\iota (x))\in \hat I_R\right\}>\eta,$$ then $\hat f^j(\iota (x))\in \hat I_R\sm B_N$ infinitely often. \label{lem:away from crit}
\end{lemma}

\begin{proof}
In a Hofbauer extension, if a point $\hat x\in \hat I$ is very close to $\bd \hat I$ then its $\hat f$-orbit shadows a point in $\bd \hat I$ for a very long time, and so it must spend a long time high up in the Hofbauer extension.  Therefore we can choose $p,N\in \N$ so that $\hat x\in B_{N}(\bd \hat I)\cap\hat I_R$ implies that \begin{equation}\hat f^{p}(\hat x)\in \hat I\sm\hat I_R \text{ and } \frac1{p}\#\{1\le j\le p:\hat f^j(\hat x)\in \hat I_R\}<\eta.\label{eq:high in tower}
\end{equation}
Suppose, for a contradiction, that $k$ is the last time that, for $x\in I$, $\hat f^k(\iota (x))\in \hat I_R\sm B_N$.  Then if $\hat f^j(\iota (x))\in \hat I_R$ for $j>k$ then $\hat f^j(\iota (x))$ must be contained in $B_N$.  Hence by \eqref{eq:high in tower}, we have
$$\limsup_k\frac1k\#\{1\le j\le k:\hat f^j(\iota (x))\in \hat I_R\}<\eta,$$ a contradiction.
\end{proof}

\begin{proof}[Proof of Proposition~\ref{prop:LE}]
We choose $R, N\in \N$, $\oG_\lambda'$ as in Lemmas~\ref{lem:big scale} and \ref{lem:away from crit} so that for any $x\in \oG_\lambda'$, $\iota (x)$ enters $\hat I_R\sm B_N$ infinitely often.

In the following we can deal with either inducing schemes of type A or type B.
We can choose $\eps>0$ so small that $\hat I_R\sm B_N\subset \cup_{\hat X\in Cover(\eps)}\hat X$.  We denote the set of points $\hat x\in \hat I$ so that the orbit of $\hat x$ enters $\hat X\subset \hat I$ infinitely often by $\hat X^\infty$.  Therefore,  for $x\in \oG_\lambda'$, there exists $\hat X_k\in Cover(\eps)$ so that $\iota (x)\in \hat X_k^\infty$.  Thus
$$\oG_\lambda'=\bigcup_{k=1}^n\{x\in \oG_\lambda':\iota (x)\in \hat X_k^\infty\}.$$
Therefore, we can choose a particular $\hat X_k$ so that
$$\dim_H(\oG_\lambda')=\dim_H\left\{x\in \oG_\lambda':\iota (x)\in \hat X_k^\infty\right\},$$
as required.
\end{proof}
\section{Proof of main results}
\label{sec:main thm}

For a potential  $\phi:I \to \R$, if the Birkhoff average $\lim_{n\to \infty}\frac{S_n\phi(x)}n$ exists, then we denote this limit by $S_\infty\negthinspace\phi(x)$.  If $\Phi$ is some induced potential, we let $S_\infty\negthinspace\Phi(x)$ be the equivalent average for the inducing scheme.

\begin{remark}  Let $f\in \F$ satisfy \eqref{eq:grow} and $\phi$ be a H\"older potential satisfying \eqref{eq:range} and $P(\phi)=0$.  Proposition~\ref{prop:eqgen} implies that there exists an equilibrium state $\mu_\phi$, but also for an inducing scheme $(X,F)$, it must have $P(\Phi)=0$ for the induced potential $\Phi$.  In fact this is only stated for type A inducing schemes in Proposition~\ref{prop:eqgen}, but will we prove this for type B schemes as well in Lemma~\ref{lem:press 0}.

For $x\in X$, we define $$\check d_{\mu_\Phi}(x):=\lim_{n\to \infty}\frac{\log\mu_{\Phi}(\c_n^F[x])}{-\log|DF^n(x)|}$$ if the limit exists.  Here $\c_n^F[x]$ is the $n$-cylinder at $x$ with respect to the inducing scheme $(X,F)$.  Since $P(\Phi)=0$, the Gibbs property of $\mu_\Phi$ implies
$$\check d_{\mu_{\Phi}}(x)= \lim_{n\to \infty}\frac{\Phi_n(x)}{-\log|DF^n(x)|}$$ whenever one of the limits on the right exists.  Also note that if both $S_\infty\negthinspace\Phi(x)$ and $\lambda_F(x)$ exist then $\check d_{\mu_\Phi}(x)$ also exists.  Suppose that $S_\infty\negthinspace\Phi(x)$ exists.  It was shown by Pollicott and Weiss \cite[Proposition 3]{PolWei} that if we also know

\begin{list}{$\bullet$}{\itemsep 0.5mm}
\setlength{\itemindent=-7mm}
%\begin{itemize}
%\begin{enumerate}
\item $\check d_{\mu_\Phi}(x)$ exists, then $d_{\mu_\Phi}(x)$ and $\lambda_F(x)$ exist and $d_{\mu_\Phi}(x)=\check d_{\mu_\Phi}(x) =\frac{S_\infty\negthinspace\Phi(x)}{-\lambda_F(x)}$;
\item $d_{\mu_\Phi}(x)$ exists, then $\check d_{\mu_\Phi}(x)$ and $\lambda_F(x)$ exist and $\check d_{\mu_\Phi}(x)= d_{\mu_\Phi}(x)=\frac{S_\infty\negthinspace\Phi(x)}{-\lambda_F(x)}$.
%\end{enumerate}
%\end{itemize}
\end{list}

Note that for $x\in (X, F)^\infty$ we can write $$\frac{\Phi_n(x)}{-\log|DF^n(x)|} = \frac{\left(\frac{\phi_{n_k}(x)}{n_k}\right)}{\left(\frac{-\log|Df^{n_k}(x)|}{n_k}\right)}$$
where $n_k=\tau^k(x)$.  Hence we can replace any assumption on the existence of $S_\infty\negthinspace\Phi(x)$ and $\lambda_F(x)$ above by the existence of $S_\infty\negthinspace\phi(x)$ and $\lambda_f(x)$.
\iffalse

As in \cite[Proposition 4]{Iom}, for the potential $\Phi$, since $P(\Phi)=0$, the Gibbs property of $m_\Phi$ implies
$$d_{\mu_{\Phi}}(x)=\lim_{n\to \infty}\frac{\log\mu_{\Phi}(\c_n[x])}{-\log|DF^n(x)|}= \lim_{n\to \infty}\frac{\Phi_n(x)}{-\log|DF^n(x)|}$$ whenever one of the limits exists.  That is to say, if $\overline{\Phi}(x)$ and $\lambda_F(x)$ exist, we have $$\overline{\Phi}(x)=-\lambda_F(x)d_{\mu_{\Phi}}(x).$$

Note that for $x\in (X, F)^\infty$ we can also write
$$d_{\mu_{\Phi}}(x)=\lim_{n\to \infty}\frac{\log\mu_{\Phi}(\c_n[x])}{-\log|DF^n(x)|}= \lim_{n\to \infty}\frac{\left(\frac{\Phi_n(x)}{\tau^n(x)}\right)}{\left( \frac{-\log|DF^n(x)|}{\tau^n(x)}\right)}= \lim_{k\to\infty} \frac{\left(\frac{\phi_{n_k}(x)}{n_k}\right)}{\left(\frac{-\log|Df^{n_k}(x)|}{n_k}\right)}$$
where $n_k=\tau^k(x)$.
So if $S_\infty\negthinspace\phi(x)$ and $\lambda_f(x)$ exist, we have $$S_\infty\negthinspace \phi\negmedspace(x)=-\lambda_f(x)d_{\mu_{\Phi}}(x).$$\fi
\label{rmk:Weiss thing}
\end{remark}

Let $$\alpha(q):=-\frac{\int\phi~d\mu_{\psi_q}}{\int\log|Df|~d\mu_{\psi_q}} = -\frac{\int\Phi~d\mu_{\Psi_q}}{\int\log|DF|~d\mu_{\Psi_q}}.$$

For the proof Theorem~\ref{thm:main spectrum} we will need two propositions relating the pointwise dimension for the induced measure and the original measure.  The reason we need to do this here is that the induced measure $\mu_\Phi$ is not, as it would be if the inducing scheme were a first return map, simply a rescaling of $\mu_\phi$.

\begin{proposition}
Given $f\in \F$ and a H\"older potential $\phi:I \to \R$ satisfying \eqref{eq:range} and $P(\phi)=0$, then there exists an equilibrium state $\mu_\phi$ and a $\phi$-conformal measure $m_\phi$ and $C_\phi>0$ so that $$\frac1{C_\phi}\le \frac{d\mu_\phi}{dm_\phi}\le C_\phi.$$
\label{prop:conf}
\end{proposition}
Notice that this implies that $d_{m_\phi}=d_{\mu_\phi}$ and, by the conformality of $m_\phi$, $d_{\mu_\phi}(x)=d_{\mu_\phi}(f^n(x))$ for all $n\in \N$.

This proposition follows from \cite{kellholder}.  However, as we mentioned in the introduction, we can also prove the existence of conformal measures under slightly different hypotheses on the map and the potential.  The class of potentials we can deal with include discontinuous potentials satisfying \eqref{eq:range}, as well as potentials $x\mapsto -t\log|Df(x)|$ for $t$ close to 1. Since this is of independent interest, we will provide a proof of this in the appendix.  A generalised version of the following result is also proved in the appendix.

\begin{proposition}
Suppose that $f\in \F$ satisfies \eqref{eq:grow} and $\phi:I \to \R$ is a H\"older potential satisfying \eqref{eq:range} and $P(\phi)=0$. For any inducing scheme $(X,F)$ either of type A or type B, with induced potential $\Phi:X\to \R$, for the equilibrium states $\mu_\phi$ for $(I,f,\phi)$ and $\mu_\Phi$ for $(X,F,\Phi)$, there exists $C_\Phi'>0$ so that $$\frac1{C_\Phi'}\le \frac{d\mu_\Phi}{d\mu_\phi}\le C_\Phi'.$$
\label{prop:induced density}
\end{proposition}

Our last step before proving Theorem~\ref{thm:main spectrum} is to show that the function $T_\phi$ as in \eqref{eq:Tphi} is strictly convex, which will mean that $\DS_\phi$ is strictly convex also, and the sets $U$ will contain non-trivial intervals.

\begin{lemma}
Suppose that $f\in \F$ satisfies \eqref{eq:grow} and $\phi$ is a H\"older potential satisfying \eqref{eq:range}. Then either there exists $\delta>0$ such that $T_\phi$ is strictly convex in $$PB(\phi)\cap\left((-\delta,\delta)\cup (1-\delta, 1+\delta)\right),$$ or $\mu_\phi=\mu_{-\log|Df|}$. \label{lem:affine}
\end{lemma}

\begin{remark}
For the particular case when $f\in \F$ and $\phi$ is a constant potential, in which case $P(\phi)=0$ implies $\phi\equiv -\htop(f)$, Lemma~\ref{lem:affine} says that $T_\phi$ is not convex if and only if $\mu_{-\log|Df|}=\mu_{-\htop(f)}$.  By  \cite[Proposition 3.1]{Dobbsvis}, this can only happen if $f$ has finite postcritical set.  We have excluded such maps from $\F$.
\label{rmk:preper}
\end{remark}

\begin{proof}[Proof of Lemma~\ref{lem:affine}]
Suppose that $T_\phi$ is not strictly convex on some interval $U$ intersecting a neighbourhood of $PB(\phi)\cap[0,1]$.  Since $T_\phi$ is necessarily convex, in $U$ it must be affine.   We will observe that for all $q\in U$, the equilibrium state for $\psi_q$ is the same.   We will then show that $[0,1]\subset U$.  Since \eqref{eq:grow} holds, and hence there is an acip $\mu_{-\log|Df|}$, this means that $\mu_\phi\equiv\mu_{-\log|Df|}$.

Our assumptions on $U$ imply that there exists $q_0\in U$ so that for a relevant inducing scheme $(X, F)$, there exists $\beta>0$ so that $\mu_{\Psi_{q_0}}\{\tau\ge n\}=O(e^{-\beta n})$.  Moreover, $DT_\phi(q)$ is some constant $\gamma\in \R$ for all $q\in U$.
As in for example \cite[Section II]{PesWei_mult} or  \cite[Chapter 7 p.211]{Pesbook} the differentiability of $T_\phi$ implies that
$\frac{\int\phi~d\mu_{\psi_q}}{\lambda(\mu_{\psi_q})}=-\gamma$ for all $q\in U$.
Since by definition $P(\psi_q)=0$, for $q_0\in U$, any measure $\mu$ with $\frac{\int\phi~d\mu}{\lambda(\mu)}=\frac{q_0}{T_\phi(q_0)}$ must be an equilibrium state for $\psi_{q_0}$.  Since there is a unique measure for $\psi_{q_0}$ we must have $\gamma=\frac{q_0}{T_\phi(q_0)}$ and $\mu_{\psi_q}=\mu_{\psi_{q_0}}$ for all $q\in U$.

By Proposition~\ref{prop:q nhd} there exists $\delta>0$ such that $(1-\delta, 1+\delta)\subset PB(\phi)$ and $(0, \delta) \subset PB(\phi)$.  If, moreover, $PB(\phi)$ contains a neighbourhood of 0 then we can adjust $\delta>0$ so that $(-\delta, \delta) \subset PB(\phi)$.

\textbf{Case 1:}
Suppose that $U \cap PB(\phi)\cap (1-\delta, 1+\delta) \neq \es$.  Since by Proposition~\ref{prop:q nhd}, $T_\phi$ is analytic in this interval, $T_\phi$ must be affine in the whole of $(1-\delta, 1+\delta)$.  Therefore $1\in U$.  We will prove that $0\in U$.  By Proposition~\ref{prop:q nhd} we can choose a type A inducing scheme $(X,F)$ so that $\mu_{\psi_q}$ is compatible with $(X,F)$ for all $q\in (1-\delta, 1+\delta)$.  Recall from Proposition~\ref{prop:eqgen} that there exists $\beta_\Phi>0$ so that $\mu_{\Psi_{1}}\{\tau\ge n\}=O(e^{-\beta_\Phi n})$.

We suppose that $0\le q<1$, and hence $T_\phi(q)\ge 0$.  We choose $q_0> 1-\delta$ very close to $1-\delta$.  Then by convexity $T_\phi(q) \ge T_\phi(q_0)+\gamma(q-q_0)$. Hence, for $Z_0^*$ as in \eqref{eq:Z0},
\begin{align*}
Z_0^*(\Psi_q)& =\sum_n n\sum_{\tau_i=n}|X_i|^{T_\phi(q)}e^{q\Phi_i} \le \sum_n n\sum_{\tau_i=n}|X_i|^{T_\phi(q_0)+\gamma(q-q_0)}e^{q\Phi_i} \\
&\le \sum_n n\sup_{\tau_i=n}\left(|X_i|^{\gamma(q-q_0)} e^{(q-q_0)\Phi_i}\right) \sum_{\tau_i=n}|X_i|^{T_\phi(q_0)}e^{q_0\Phi_i} \\
&\le \sum_nn e^{n(q-q_0)\inf\phi} \sum_{\tau_i=n}|X_i|^{T_\phi(q_0)}e^{q_0\Phi_i}.\label{eq:linear}\end{align*}
By the Gibbs property of $\mu_{\Psi_{q_0}}$, we can estimate $\sum_{\tau_i=n}|X_i|^{T_\phi(q_0)}e^{q_0\Phi_i}$ by $\mu_{\Psi_{q_0}}\{\tau=n\}=\mu_{\Psi_{1}}\{\tau=n\} \le e^{-\beta_\Phi n}$.  So if $(q-q_0)\inf\phi<\beta_\Phi$ then similarly to the proof of Proposition~\ref{prop:q nhd}, $q\in PB(\phi)$.  Since $T_\phi$ is analytic in $PB(\phi)$, this means that $T_\phi$ is still affine at $q$ and therefore that $U$ was not the largest domain of affinity `to the left'.  We can continue doing this until we hit the left-hand boundary of $PB(\phi)$.  In particular, this means that  $0\in U$.

\textbf{Case 2:}
Suppose that $PB(\phi)\cap(-\delta,\delta)\cap U\neq \es$.  As in Case 1, this implies $[0, \delta]\subset U$.  We will prove that $1\in U$.

By Proposition~\ref{prop:q nhd} we can choose a type B inducing scheme $(X,F)$ so that $\mu_{\psi_q}$ is compatible with $(X,F)$ for all $q\in (\delta', \delta)$ where $\delta':=\delta/2$.  Recall from Proposition~\ref{prop:eqgen} that there exists $\beta_{DF}>0$ so that $\mu_{\Psi_{\delta'}}\{\tau\ge n\}=O(e^{-n\beta_{DF}})$.

We let $\delta<q\le 1$ and $q_0<\delta$ be very close to $\delta$.
Again by convexity $T_\phi(q) \ge T_\phi(q_0)+\gamma(q-q_0)$.  Similarly to Case 1, \begin{align*}
Z_0^*(\Psi_q)& =\sum_nn\sum_{\tau_i=n}|X_i|^{T_\phi(q)}e^{q\Phi_i} \le \sum_nn\sum_{\tau_i=n}|X_i|^{T_\phi(q_0)+\gamma(q-q_0)}e^{q\Phi_i} \\
&\le \sum_nn\sup_{\tau_i=n}\left(|X_i|^{\gamma(q-q_0)} e^{(q-q_0)\Phi_i}\right) \sum_{\tau_i=n}|X_i|^{T_\phi(q_0)}e^{q_0\Phi_i}.\end{align*}
Since $|X_i|\ge e^{-\tau_i|Df|_\infty}$,
$$\sup_{\tau_i=n}\left(|X_i|^{\gamma(q-q_0)} e^{(q-q_0)\Phi_i}\right)\le e^{n(q-q_0)(-\gamma |Df|_{\sup}+\sup\phi) }.$$
So if $(q-q_0)(-\gamma |Df|_{\infty}+\sup\phi)<\beta_{DF}$ then similarly to Case 1 we can  conclude that all points in $PB(\phi)$ to the right of $q_0$ are in $U$.  In particular $1\in U$.

In both cases 1 and 2, we concluded that $[0,1]\subset U$.  Therefore $\mu_\phi\equiv\mu_{-\log|Df|}$, as required.
\end{proof}

\begin{proof}[Proof of Theorem~\ref{thm:main spectrum}] Let $L_\phi$ be the Legendre transform of $T_\phi$ as in \eqref{eq:Tphi} wherever these functions are defined.

\textbf{The upper bound: \boldmath $\widetilde{\DS}_\phi\le L_\phi$\unboldmath.}
To get this bound, we first pick a suitable inducing scheme.  Given $q\in PB(\phi)$, since $\tilde \K_\phi(\alpha(q))=\cup_{n\ge 1} \oG_{\frac1n}\cap \tilde \K_\phi(\alpha(q))$, for all $\eta>0$ there exists $\lambda>0$ so that $\dim_H(\oG_\lambda'\cap \tilde \K_\phi(\alpha(q)))\ge \dim_H(\tilde \K_\phi(\alpha(q)))-\eta$.  For some $s<\dim_H(\tilde \K_\phi(\alpha(q)))$, we take an inducing scheme $(X,F)$ as in Proposition~\ref{prop:LE} (this can be for schemes of type A or B, whichever we need).

We next show that $\widetilde{\DS}_\phi\le \DS_\Phi$ and then use Theorem~\ref{thm:Iommi} and Proposition~\ref{prop:q nhd} to conclude the proof of the bound.
Let $x\in \K_\phi(\alpha)\cap \oG_\lambda'$.  By transitivity there exists $j$ so that $x\in f^j(X)$.  Let $y\in X$ be such that $f^j(y)=x$.  Since $x\in \oG_\lambda'$, we must also have $y\in (X,F)^\infty$ by Proposition~\ref{prop:LE}.  By Propositions~\ref{prop:conf} and \ref{prop:induced density}, $d_{\mu_\phi}(x)=d_{\mu_\phi}(y)=d_{\mu_\Phi}(y)$, so $y\in  \K_\Phi(\alpha)$.  Therefore,
$$\tilde \K_\phi(\alpha)\cap \oG_\lambda'\subset \cup_{k=0}^\infty f^k(\K_\Phi(\alpha)).$$
Hence $$\widetilde{\DS}_\phi-\eta \le \dim_H(\K_\phi(\alpha)\cap \oG_\lambda') \le \dim_H\left(\cup_{k=0}^\infty f^k(\K_\Phi(\alpha))\right).$$  Since $f$ is clearly Lipschitz, $\dim_H\left(\cup_{k=0}^\infty f^k(\K_\Phi(\alpha))\right)= \dim_H(\K_\Phi(\alpha))$, so $\widetilde{\DS}_\phi(\alpha)-\eta \le \DS_\Phi(\alpha)$.  Theorem~\ref{thm:Iommi} says that $\DS_\Phi(\alpha(q))$ is $L_\Phi(\alpha)$, the Legendre transform of $T_\Phi$.  Therefore, $\widetilde{\DS}_\phi(\alpha)-\eta \le L_\Phi(\alpha)=L_\phi(\alpha)$, where the final equality follows from Proposition~\ref{prop:q nhd}.  Since $\eta>0$ was arbitrary, we have $\widetilde{\DS}_\phi(\alpha)\le L_\phi(\alpha)$.

\textbf{The lower bound: \boldmath $\widetilde{\DS}_\phi \ge L_\phi$\unboldmath.}
We will use the Hausdorff dimension of the equilibrium states for $\psi_q$ to give us the required upper bound here.
For $\mu\in \M_+$, by Theorem~\ref{thm:schemes} there exists an inducing scheme $(X,F)$ which $\mu$ is compatible to.  This can chosen to be of type A or type B. By Proposition~\ref{prop:induced density}, $d_{\mu_\phi}(x)=d_{\mu_\Phi}(x)$ for any $x\in (X,F)^\infty$, where $\Phi$ is the induced potential for $(X,F)$.
Now suppose that $\frac{\int\phi~d\mu}{\lambda_f(\mu)}=-\alpha$.  Then for $\mu$-\ae $x$, $S_\infty\negthinspace\phi(x)$ and $\lambda(x)$ exist, and by the above and Remark~\ref{rmk:Weiss thing}, since we may choose $X$ so that for $x\in (X, F)^\infty$, we have $$d_{\mu_\phi}(x)=d_{\mu_\Phi}(x)= \frac{S_\infty\negthinspace\phi(x)}{-\lambda_f(x)} =\alpha.$$  Hence $\mu$-\ae $x$ is in
$\K_\phi(\alpha)$.
Therefore,
$$\widetilde{\DS}_\phi(\alpha)\ge \sup\left\{\frac{h_\mu}{\lambda_f(\mu)}:\mu\in \M_+ \text{ and } \frac{\int\phi~d\mu}{\lambda_f(\mu)}=-\alpha\right\}.$$
By Lemma~\ref{lem:PB Scover}, we know that there is an equilibrium state $\mu_{\psi_q}$ for $\psi_q$.  Then by definition, $h_{\mu_{\psi_q}}+\int-T(q)\log|Df|+q\phi~d\mu_{\psi_q}=0$.  Therefore, for $\alpha=\alpha(q)$,
$$\frac{h_{\mu_{\psi_q}}}{\lambda_f(\mu_{\psi_q})}=T(q)+q\alpha=L_\phi(\alpha).$$
And hence $\widetilde{\DS}_\phi(\alpha)\ge L_\phi(\alpha)$.
Putting our two bounds together, we conclude that $\widetilde{\DS}_\phi(\alpha)=L_\phi(\alpha)$.

We next show (a), (b) and (c).  First note that since we have assumed that $\mu_\phi\neq\mu_{-\log|Df|}$, Lemma~\ref{lem:affine} means that $T_\phi$ is strictly convex in $PB(\phi)$.  This implies that $U$ will contain non-trivial intervals.  For example, if \eqref{eq:grow} holds then $P(\phi)=0$ and \cite{HofRaidim} imply that
$$\alpha(1)= -\frac{\int\phi~d\mu_{\phi}}{\lambda_f(\mu_{\phi})}= \frac{h_{\mu_{\phi}}}{\lambda_f(\mu_{\phi})}= \dim_H(\mu_\phi).$$
By Proposition~\ref{prop:q nhd} and Lemma~\ref{lem:affine}, for any $\alpha$ close to $\dim_H(\mu_\phi)$ there exists $q$ \st $DT_\phi(q)=\alpha$.  Hence by the above,
$\widetilde{\DS}_\phi(\alpha)=L_\phi(\alpha)$.

Similarly, let us assume that \eqref{eq:poly} holds.  We have $$\alpha(0) = -\frac{\int\phi~d\mu_{-\log|Df|}}{\lambda_f(\mu_{-\log|Df|})}=\alpha_{ac}.$$ So the arguments above, Proposition~\ref{prop:q nhd} and Lemma~\ref{lem:affine} imply that for any $\alpha<\alpha_{ac}$ there exists $q$ \st $DT_\phi(q)=\alpha$, and also
$\widetilde{\DS}_\phi(\alpha)=L_\phi(\alpha)$.
The same holds for all $\alpha$ in a neighbourhood of $\alpha_{ac}$ when \eqref{eq:CE} holds.
\end{proof}

\begin{proof}[Proof of Proposition~\ref{prop:nonreg}.]
It was pointed out in \cite[Remark 4.9]{Iom} that by \cite{BarSch}, for an inducing scheme $(X,F)$ with potential $\Phi:X\to \R$, the Hausdorff dimension of the set of points with $d_{\mu_\Phi}(x)$ not defined has the same dimension as the set of points for which the inducing scheme is defined for all time.  So we can choose $(X,F)$ to be any inducing scheme which is compatible to the acip to show that the Hausdorff dimension of this set of points is 1.  In fact any type A or type B inducing scheme is compatible to the acip.  By Proposition~\ref{prop:induced density}, if $d_{\mu_\Phi}(x)$ not defined then neither is $d_{\mu_\phi}(x)$, so the proposition is proved.
\end{proof}

\subsection{Going to large scale: the proof of Corollary~\ref{cor:poly spectrum}}
\label{sec:GrSm}

Suppose that $f\in \F$ extends to a polynomial on $\mathbb{C}$ with no parabolic points and all critical points in $I$.  In the context of rational maps, Graczyk and Smirnov \cite{GrSm} prove numerous results for such maps satisfying \eqref{eq:poly}.  For $\delta>0$, we say that $x$ \emph{goes to $\delta$-large scale at time $n$} if there exists a neighbourhood $W$ of $x$ such that $f:W \to B_\delta(f^n(x))$ is a diffeomorphism.  It is proved in \cite{GrSm} that there exists $\delta>0$ such that the set of points which do not go to $\delta$-large scale for an infinite sequence of times has Hausdorff dimension less than $\frac{\ell_{max}(f)}{\betaP-1}<1$ where $\betaP$ is defined in \eqref{eq:poly}.  Here we will sketch how this implies Corollary~\ref{cor:poly spectrum}.

By \cite{Kellift}, if $f\in \F$ and $x\in I$ goes to $\delta$-large scale with frequency $\gamma$, then there exists $N=N(\delta)$ so that iterates of $\iota (x)$ by $\hat f$ enter $\hat I_N$ with frequency at least $\gamma$.  In \cite{Kellift, BTholo}, this idea was used to prove that for $\mu\in \M_{erg}$, if $\mu$-\ae $x$ goes to $\delta$-large scale with some frequency greater than $\gamma>0$, then there exists $\hat\mu$ an ergodic $\hat f$-invariant probability measure on $\hat I$, with $\hat\mu(\hat I_N)>\gamma$ (so also $\hat\mu$-\ae $\hat x$ enters $\hat I_N$ with positive frequency), and $\mu=\hat\mu\circ\pi^{-1}$.  By the arguments above this means that we can build an inducing $(X,F)$ scheme from a set $\hat X\in \hat I_N$ which is compatible to $\mu$.

However, to prove Corollary~\ref{cor:poly spectrum}, we only need that sufficiently many points $x$ have $k\ge 0$ such that $f^k(x)\in (X,F)^\infty$, which does not necessarily mean that these points must go to large scale with positive frequency.  (Note that we already know that all the measures $\mu$ we are interested in can be lifted to $\hat I$.) We only need to use the fact, as above, that if $A$ is the set of points which go to $\delta$-large scale infinitely often, then there exists $R\in \N$ so that for all $x\in A$, $\iota (x)$ enters $\hat I_R$ infinitely often.  Hence the machinery developed above `sees' all of $A$, up to a set of Hausdorff dimension $<\frac{\ell_{max}(f)}{\beta_P-1}$.  Since this value is $<1$, for our class of rational maps, we have $\DS_\phi(\alpha)=\widetilde{\DS}_\phi(\alpha)$ for $\alpha$ close to $\alpha_{ac}$.  Similarly, if $\frac{\ell_{max}(f)}{\beta_P-1}<\dim_H(\mu_\phi)$ then the same applies for $\alpha$ close to $\dim_H(\mu_\phi)$.

Note that for rational maps as above, but satisfying \eqref{eq:CE}, the same argument gives another proof of Corollary~\ref{cor:CE spectrum}.

It seems likely that the analyticity condition can be weakened to include all maps in $\F$ satisfying \eqref{eq:poly}.

\subsection{Points with zero Lyapunov exponent can be seen}
\label{sec:what is seen}

In this section we discuss further which points can and cannot be seen by the inducing schemes we use here.

Suppose that $(X, F, \tau)$ is an inducing scheme of type A.  Then there is a corresponding set $\hat X\subset \hat I$ such that $\tau(y)$ is $r_{\hat X}(\hat y)$ where $\hat y\in \hat X$ is such that $\pi(\hat y)=y$ and $r_{\hat X}$ is the first return time to $\hat X$.  Then there exist points $\hat x\in \hat X$ so that $\pi(\hat f^k(x))\in \crit$ and $\hat f^j(\hat x)\notin \hat X$ for all $1\le j< k$.  This implies that from iterate $k$ onwards, this orbit is always in the boundary of its domain $D\in \D$.  Since $\hat X$ is always chosen to be compactly contained inside its domain $D_{\hat X}\in \D$, this means that $\hat x$ never returns to $\hat X$.  Hence for $x=\pi(\hat x)$, $\tau(x)=\infty$.  On the other hand, there are precritical points $x$ with $\hat x=\pi|_{\hat X}^{-1}(x)$ which returns to $\hat X$ before it hits a `critical line' $\pi^{-1}(c)$ for $c\in \crit$.  For such a point, $\tau(x)<\infty$, but for all large iterates $k$, we must have $\tau(f^k(x))=\infty$.  Hence precritical points in $X$ cannot have finite inducing time for all iterates.  This can be shown similarly for type B inducing schemes.
\iffalse Note also that since $X$ was chosen to be a cylinder set, \ie to be in some $\P_n$, then any precritical point  $x$ is a precritical point and $f^k(x)$ is the first iterate $\ge 0$ for which $f^k(x)\in X$, then $k\ge n$\fi  We can extend this to show that no precritical point is counted in our proof of Theorem~\ref{thm:main spectrum}.

Moreover, in this paper we are able to find $\widetilde{\DS}_\phi(\alpha)$ through measures on $\K_\alpha$.  In fact we can only properly deal with measures which are compatible to some inducing scheme.  As in Theorem~\ref{thm:schemes}, the only measures we can consider are in $\M_+$.  This means that the set of points $x$ with $\underline\lambda(x)=0$ is not seen by these measures.  As pointed out above Corollary~\ref{cor:CE spectrum}, \cite{BrvS} shows that in the Collet-Eckmann case, the set of points with $\overline\lambda(x)=0$ is countable and thus has zero Hausdorff dimension.  (Note that even in this well-behaved case it is not yet clear that the set of points with $\underline\lambda(x)=0$ has zero Hausdorff dimension.)  The general question of what is the Hausdorff dimension of $I\sm \oG_0$ for topologically transitive maps is, to our knowledge, open.

On the other hand, it is not always the case that given an inducing scheme $(X, F, \tau)$, all points $x\in X$ for which $\tau(F^k(x))<\infty$ for all $k\ge 0$ have positive Lyapunov exponent.  For example, we say that $f$ has \emph{uniform hyperbolic structure} if $\inf\{\lambda_f(p):p \text{ is periodic}\}>0$.  Nowicki and Sands \cite{NoSa} showed that for unimodal maps in $\F$ this condition is equivalent to \eqref{eq:CE}.  If we take $f\in \F$ without uniform hyperbolic structure, then it can be shown that for any inducing scheme $(X,F,\tau)$ as above, there is a sequence $\{n_k\}_{k}$ such that $$\frac{\sup\{\log|DF(x)|:x\in X_{n_k}\}}{\tau_{n_k}}\to 0.$$
There exists $x\in X$ so that $F^k(x)\in X_{n_k}$ for all $k$.  Thus $\underline\lambda(x)\le 0$, but $\tau(F^k(x))<\infty$ for all $k\ge 0$.  In the light of the proof of Corollary~\ref{cor:poly spectrum}, we note that $x$ goes to $|X|$-large scale infinitely often, but with zero frequency.

\iffalse
suppose that there exists $C>0$ such that $\frac{n^{-\beta}}C\le \tau_{-\log|DF|}\{\tau=n\} \le C n^{-\beta}$ for some $\beta>1$.  Since not much is currently known for lower bounds on such tails, one could argue that this is an artificial situation.  However, it should be expected that at least in the unimodal case, maps satisfying \eqref{eq:poly} have such behaviour. For convenience we can label the inducing domains $X_i$ in order of increasing inducing time, \ie $i'>i$ implies $\tau_{i'}\ge \tau_i$.  Then there exists $x\in X$ so that $F^i(x)\in X_i$ for all $i$.  Thus, since $|DF_i(F^{i-1}(x))| \asymp \frac{|X|}{|X_i|} =O(\tau_i^\beta)$, $\frac{\log|DF^n(x)|}{\tau^n} \to 0$.  So $\underline\lambda(x)\le 0$, but $\tau(F^n(x))<\infty$ for all $n\ge 0$.  In the light of the proof of Corollary~\ref{cor:poly spectrum}, we note that $x$ goes to $|X|$-large scale infinitely often, but with zero frequency.
\fi

In conclusion, while it may not be necessary, it seems to be extremely difficult to study notions such as dimension spectra unless we are allowed to exclude points $x$ with $\overline\lambda(x)\le 0$ from consideration.

\iffalse
\textbf{Proof of the above bit}

First let's deal only with $\hat f$-periodic points.  Suppose that $\hat X$ is the domain we're taking to give us our inducing scheme.  Suppose that we do not have a sequence of $\hat f$-periodic points $\{\hat p_n\}_n\subset \hat X$ so that $\lambda_{\hat f}(\hat p_n) \to 0$.  This means that we have to look outside $\hat X$ for such points.  Take a periodic point $\hat p\in \hat I_{\T}$ of period $q$ and a small neighbourhood $U\ni \hat p$.  Then by transitivity, there exists $n_1, n_2\in \N$ and $\hat x\in \hat X$ so that $\hat f^{n_1}(X) \supset U$, $\hat f^{n_1}(\hat x)\in U$ and $n_1+n_2$ is the first return time of $\hat x$ to $\hat X$.  We can adjust $U$ so that $n_1$ is the first time $\hat f^k(\hat X)\cap U\neq \es$.  Then since $\hat p$ is a repelling periodic point with orbit disjoint from $\hat X$, for all $k\in \N$ there is $\hat y_k$ so that $\hat f^{kq}(\hat y_k)=\hat f^{n_1}(\hat x)$, $\hat f^{jq}(\hat y_k) \in U$ for all $0\le j\le k$ and $\hat f^{j}(\hat y_k) \notin \hat X$ for all $0\le j\le kq$.  By the setup, for all such $k$, there exists $\hat x_k$ such that $\hat f^{n_1}(\hat x_k)=y_k$.  Around each such $x_k$ there will be a domain $\hat X_{i_k}$.  Choosing $U$ small (so that $D\hat f^q|_U\approx \lambda_{\hat f}(\hat p)$), for $k$ large enough, we have $\frac{\log |DF|_{X_{i_k}}}{\tau_{i_k}}$ arbitrarily close to $\lambda_{\hat f}(\hat p)$, which can in turn be taken arbitrarily close to 0.

\fi

\section{Lyapunov spectrum}
\label{sec:LE spec}

For $\lambda\ge 0$ we let
$$L_\lambda=L_\lambda(f):=\left\{x: \lambda_f(x)=\lambda \right\} \text{ and } L'=L'(f):=\left\{x: \lambda_f(x) \text{ does not exist}\right\}.$$
The function $\lambda\mapsto \dim_H(L_\lambda)$ is called the \emph{Lyapunov spectrum}.
Notice that by \cite{BrvS}, if $f\in \F$ satisfies \eqref{eq:grow} then if the Lyapunov exponent at a given point exists then it must be greater than or equal to 0.
In this section we explain how the results above for pointwise dimension are naturally related to the Lyapunov spectrum.  As we show below, the equilibrium states $\mu_{-t\log|Df|}$ found in \cite{PeSe, BTeqnat} for certain values of $t$, depending on the properties of $f$, are the measures of maximal dimension sitting on the sets $L_\lambda$ for some $\lambda=\lambda(t)$.

Recall that $\mu_{-\log|Df|}$ is the acip for $f$.  We denote the measure of maximal entropy by $\mu_{-\htop(f)}$ since it is the equilibrium state for a constant potential $\phi_a(x)=a$ for all $x\in I$; and in order to ensure $P(\phi_a)=0$, we can set $a=-\htop(f)$.  We let $\DS_{-\htop(f)}(\alpha)=\dim_H(\K_{-\htop(f)}(\alpha))$ where $\K_{-\htop(f)}$ is defined for the measure $\mu_{-\htop(f)}$ as above.

\begin{proposition}
If $f\in \F$ then there exists an open set $U\subset \R$ containing $\frac{\htop(f)}{\lambda_f(\mu_{-\htop(f)})}$ so that for each  $\alpha\in U$ the values of $\dim_H\left(L_{\frac{\htop(f)}\alpha}\right)=\DS_{-\htop(f)}(\alpha)$
are given as the Legendre transform of $T_{-\htop(f)}$ at $\alpha$.  If $f$ satisfies \eqref{eq:poly}, then $\frac{\htop(f)}{\lambda_f(\mu_{-\log|Df|})}$ is in the closure of $U$, and if $f$ satisfies \eqref{eq:CE} then $\frac{\htop(f)}{\lambda_f(\mu_{-\log|Df|})}$ is contained in $U$.
\label{prop:lyap spec}
\end{proposition}

As observed by Bohr and Rand, this proposition would have to be adapted slightly when we are dealing with quadratic Chebyshev polynomial (which is not in our class $\F$).  In this case, $\mu_{-\htop(f)}=\mu_{-\log|Df|}$, so the Lyapunov spectrum can not analytic in a neighbourhood of 1.  Note that this agrees with Lemma~\ref{lem:affine} and Remark~\ref{rmk:preper}.

Note that the first part of the proposition makes no assumption on the growth of $|Df^n(f(c))|$ for $c\in \crit$.
The proof of this proposition follows almost exactly as in the proof of Proposition~\ref{prop:q nhd}, so we only give a sketch.

\begin{proof}
Given an inducing scheme $(X,F)$, by Remark~\ref{rmk:Weiss thing}, for each $x\in (X, F)^\infty$ if $\lambda_f(x)$ exists then
$$\lambda_f(x)=\frac{\htop(f)}{d_{\mu_{-\tau\htop(f)}}(x)}.$$  Here the potential is $\phi\equiv-\htop(f)$, and the induced potential is $-\tau\htop(f)$. This means that we can get the Lyapunov spectrum directly from $d_{\mu_{-\tau\htop(f)}}$.  As in Proposition~\ref{prop:induced density}, $d_{\mu_{-\tau\htop(f)}}(x)=d_{\mu_{-\htop(f)}}(x)$ for all $x\in X$.

Therefore it only remains to discuss the interval $U$, i.e. the equivalent of Proposition~\ref{prop:q nhd}.  First we note that Lemma~\ref{lem:affine} holds in this case without any assumption on the proof of $|Df^n(f(c))|$ for $c\in \crit$.
We fix an inducing scheme $(X,F)$.  That $Z_0^*(\Psi_q+\delta_q\tau)<\infty$ for some small $\delta_q>0$, for $q$ in some open interval $U$ can be proved exactly in the same way as in the proof of Proposition~\ref{prop:q nhd}.
\end{proof}

Note that similarly to Proposition~\ref{prop:nonreg}, the set of points for which the Lyapunov exponent is not defined has Hausdorff dimension 1.

\begin{remark}
For $t\in \R$, let $P_t:=P(-t\log|Df|)$.  It follows that $P_{T_{-\htop(f)}(q)}=q\htop(f)$.  Since $\mu_{\psi_q}$ is an equilibrium state for $-T_{-\htop(f)}(q)\log|Df|-q\htop(f)$, then it is also an equilibrium state for $-T_{-\htop(f)}(q)\log|Df|$.  Therefore, the measures for $\psi_q$ are precisely those found for the potential $-t\log|Df|$ in Proposition~\ref{prop:eqnat} and in \cite[Theorem 6]{BTeqgen}.
\end{remark}

\begin{remark}
If \eqref{eq:CE} does not hold, then Proposition~\ref{prop:lyap spec} does not deal with $L_\lambda$ for $\lambda< \lambda(\mu_{-\log|Df|})$.  This is because, at least in the unimodal case, we have no equilibrium state with positive Lyapunov exponent for the potential $x\mapsto -t\log|Df(x)|$ for $t> 1$ (\ie there is a phase transition at 1).

Nakaishi \cite{Nak} and Gelfert and Rams \cite{GelRams} consider the Lyapunov spectrum for Manneville-Pomeau maps with an absolutely continuous invariant measure, which has polynomial decay of correlations.  Despite there being a phase transition for $t\mapsto P_t$ at $t=1$, they are still able to compute the Lyapunov spectrum in the regime $\lambda\in [0,\lambda(\mu_{-\log|Df|}))$.  Indeed they show that $\dim_H(L_\lambda)=1$ for all these values of $\lambda$.  In forthcoming work we will show that we have the same phenomenon in our setting when \eqref{eq:poly}, but not \eqref{eq:CE}, holds.
\label{rmk:phase trans}
\end{remark}

\begin{remark}
If \eqref{eq:CE} holds then it can be computed that in the above proof, $Z_0^*(\Psi_q+\delta\tau)<\infty$ whenever $(1-T_{-\htop(f)}(q)-q)\htop(f) - \alpha T_{-\htop(f)}(q)$, where $\alpha$ is the rate of decay of $\mu_{-\log |DF|}\{\tau>n\}$ and $\delta$ is some constant $>0$. If $f$ is a Collet-Eckmann map very close to the Chebyshev polynomial, then $t\mapsto P(-t\log|Df|)$ is close to an affine map, and thus $T_{-\htop(f)}$ is also close to an affine map, then $Z_0^*(\Psi_q+\delta_q\tau)<\infty$ for all $q$ in a neighbourhood of $[0,1]$ and for some $\delta_q>0$.

The unimodal maps considered by Pesin and Senti \cite{PeSe} have the above property and so there exists $\eps>0$ so that $[0,1]\subset PB_\eps(-\htop(f))$.  However, this may not be the whole spectrum.
\end{remark}

In \cite{PeSe}, they ask if it is possible to find a unimodal map $f:I\to I$ so that there is a equilibrium state for the potential $x\mapsto -t\log|Df|$ for all $t\in (-\infty, \infty)$, and that the pressure function $t\mapsto P(-t\log|Df|)$ is analytic in this interval.  This would be in order to implement a complete study of the thermodynamic formalism.  As Dobbs points out in \cite{Dobbsphase}, in order to show this, even in the `most hyperbolic' cases, one must restrict attention to measures on a subset of the phase space: otherwise we would at least expect a phase transition in the negative spectrum.

\section*{Appendix}
\label{sec:density}

In this appendix we introduce a class of potentials for which the results in the rest of the paper hold.  We will also prove slightly generalised versions of Propositions~\ref{prop:conf} and \ref{prop:induced density}.

Given a potential $\phi$, and an inducing scheme $(X,F)$ of type A or B, as usual we let $\Phi$ be the induced potential.
If
\begin{equation}\label{SVI}
\sum_nV_n(\Phi) < \infty,
\end{equation}
then we say that $\phi$ satisfies the summable variations for induced potential condition, with respect to this inducing scheme.  If $\phi$ satisfies this condition for every type A or B inducing scheme $(X,F)$ with $|X|$ sufficiently small, we write $\phi\in SVI$.  Note that in \cite[Lemma 3]{BTeqgen} it is proved that if $\phi$ is H\"older and $f\in \F$ satisfies \eqref{eq:range} then $\phi\in SVI$.  Also in \cite{BTeqgen} it was proved that Proposition~\ref{prop:eqgen} holds for all potentials in $SVI$ satisfying \eqref{eq:range}, with no assumptions on the growth along the critical orbits.

Proposition~\ref{prop:conf} is already known in the case that $\phi$ is H\"older.  For interest, we will change the class of potentials in that proposition to those in $SVI$ satisfying \eqref{eq:range}, as well as to potentials of the form $x\mapsto -t\log|Df(x)|$.  We also widen the class of potentials considered in Proposition~\ref{prop:induced density}.  We will refer to Propositions~\ref{prop:conf} and \ref{prop:induced density}, but with only the assumptions that $f\in \F$ and $\phi\in SVI$, as Propositions~\ref{prop:conf}' and \ref{prop:induced density}'.  Note that Proposition~\ref{prop:induced density}' plus \cite[Lemma 3]{BTeqgen} implies Proposition~\ref{prop:induced density}.  The proof of these propositions requires three steps:
\begin{list}{$\bullet$}{\itemsep 0.2mm \topsep 0.2mm \itemindent -0mm \leftmargin=5mm}
\item Proving the existence of a conformal measure $m_\phi$ for a potential $\phi\in SVI$ satisfying \eqref{eq:range} and $P(\phi)=0$.  Since we do this using the measure $m_\Phi$ from Proposition~\ref{prop:eqgen}, we only really need to prove this for inducing schemes of type A.  However, it is of independent interest that this step can also be done for the potential $x\mapsto -t\log|Df(x)|-P(-t\log|Df|)$, so we allow type B inducing schemes also.
\item Proving that a rescaling of the measure $m_\phi$ is also conformal for our inducing schemes.  This will be used directly in the proof of Proposition~\ref{prop:conf}', so must hold for both type A and type B inducing schemes.  Note that this step works for all of the types of potential mentioned above.
\item Proving that the density $\frac{d\mu_\phi}{dm_\phi}$ is bounded.  We will use type A inducing schemes to prove this.   In this step, we must assume that $\phi$ is in $SVI$, satisfies \eqref{eq:range} and $P(\phi)=0$.
\end{list}

The necessary parts of the first and third of these steps are the content of Proposition~\ref{prop:conf}'.
\iffalse
the main step is to prove the existence of conformal measures $m_\phi$; then to prove that this conformal measure is also conformal for all inducing schemes; and finally to show the bounded density.  The first and final part only really requires type A inducing schemes, while the second must work for both type A and type B schemes.  Note that our proof of the existence of conformal measure also extends to potentials other than those satisfying \eqref{eq:range}, for example for the potentials $x\mapsto -t\log|Df(x)|-P(-t\log|Df|)$ as in Proposition~\ref{prop:eqnat}.
\fi
As above, for the proof  of this proposition, we only need to use type A inducing schemes.  But we will give the proof of the existence of the conformal measure for both types of schemes for interest.
%For the  density bound, we will only consider type A schemes.
Our inducing scheme $(X,F,\tau)$ is derived from a first return map to a set $\hat X\subset \hat I$.  Recall that if we have a type A scheme, then $\hat X$ is an interval in a single domain $\hat X\subset D\in \D$ in the Hofbauer extension.  In the type B case, $\hat X$ may consist of infinitely many such intervals.  We let $r_{\hat X}$ be the first return time to $\hat X$ and $R_{\hat X}=\hat f^{r_{\hat X}}$.  We let $\hat X_i$ denote the first return domains of $R_{\hat X}$.

We let $\hat\phi:=\phi\circ\pi$, and $\hat\mu_{\phi, \hat X}:=\frac{\hat\mu_\phi|_{\hat X}}{\hat\mu(\hat X)}$ be the conditional measure on $\hat X$.  As explained in \cite{BTeqnat}, the measure $\mu_\Phi$ is the same as $\hat\mu_{\phi,\hat X}\circ\pi^{-1}$.
Proposition~\ref{prop:eqgen} implies that for type A inducing schemes $(X,F)$, the induced potential $\Phi$ has $P(\Phi)=0$, and there a conformal measure and equilibrium state $m_\Phi$ and $\mu_\Phi$ and $C_\Phi>0$ so that $\frac1{C_\Phi} \le \frac{d\mu_\Phi}{dm_\Phi}\le C_\Phi$.
We show in Lemma~\ref{lem:press 0} that this is also true for type B inducing schemes.

We define $\hat m_\phi|_{\hat X}:=m_{\Phi}\circ\pi|_{\hat X}$.  We can propagate this measure throughout $\hat I$ as follows.

\iffalse
Given $\hat A\subset \hat f^k(\hat X_i)$ for $0\le k\le r_{\hat X}|_{\hat X_i}-1$, we let $\hat A'$ be defined so that $\hat f^k(\hat A')=\hat A$.  We define $\hat m_\phi(\hat A):=\int_{\hat A'}e^{-\hat\phi_k}~d\hat m_\phi$.
\fi

For $\hat x\in \hat X$ with $r_{\hat X}(\hat x)<\infty$, for $0\le k\le r_{\hat X}(\hat x)-1$, we define $$d\hat m_\phi(\hat f^k(\hat x))=e^{-\hat\phi_k(\hat x)}d\hat m_\phi|_{\hat X}(\hat x).$$

Let $(X,f)$ be a dynamical system and $\phi:X\to \R$ be a potential. We say that a measure $m$, is \emph{$\phi$-sigma-conformal} for $(X,f)$ if for any Borel set $A$ so that $f:A \to f(A)$ is a bijection,
$$m(f(A))=\int_Ae^{-\phi}~dm.$$  Or equivalently $dm(f(x))=e^{-\phi(x)} dm(x)$.
So the usual conformal measures are also sigma-conformal, but this definition allows us to deal with infinite measures.  The next two lemmas apply to potentials $\phi\in SVI$ satisfying \eqref{eq:range} and $P(\phi)=0$, or of the form $x\mapsto -t\log|Df(x)|-P(-t\log|Df|)$ as in Proposition~\ref{prop:eqnat}.

\begin{lemma}  Suppose that $(X,F)$ is a type A or type B system and $P(\Phi)=0$.
\begin{itemize}
\item[(a)]
$\hat m_\phi$ as defined above is a $\phi$-sigma-conformal measure.
\item[(b)]
Given a $\hat\phi$-sigma-conformal measure $\hat m_\phi'$ for $(\hat I, \hat f)$, then up to a rescaling, $\hat m_\phi'=\hat m_\phi$.
\end{itemize}
\label{lem:sig conf}
\end{lemma}

\begin{proof}
We first prove (a).
The $\Phi$-conformality of $m_\Phi$ implies that $\hat m_\phi|_{\hat X}$ is $\hat\Phi$-conformal for the system $(\hat X,R_{\hat X},\hat\Phi)$ for $\hat\Phi(\hat x):=\Phi(\pi(\hat x))$.

Given $\hat x\in \hat X$, if $0 \le j < r_{\hat X}(\hat x) - 1$, then the
relation $$d\hat m_{\phi} \circ \hat f(\hat f^j(\hat x)) = e^{-\hat\phi(\hat x)}
d\hat m_{\phi}(\hat f^j(\hat x))$$ is immediate from the definition. For $j = r_{\hat X}(\hat x)-1$, then $\hat f(\hat f^j(\hat x)) = R_{\hat X}(\hat x)$ and we obtain, for $\hat x\in \hat X$,
\begin{align*}
d\hat m_{\phi} \circ \hat f(\hat f^j(\hat x)) &= e^{-\hat\phi_j(\hat x)}d\hat m_{\phi}(\hat x) = d\hat m_\phi(\hat R(\hat x)) = e^{-\hat\Phi(\hat x)} d\hat m_\phi(\hat x) \\
&= e^{-\hat\phi(\hat f^{\hat r_{\hat X}(\hat x)-1}(\hat x))} e^{-\hat\phi_{r_{\hat X}(\hat x)-2}(\hat x)} d\hat m_\phi(\hat x) \\
&= e^{-\hat\phi(\hat f^{\hat r_{\hat X}(\hat x)-1}(\hat x))} d\hat m_{\phi}(\hat f^{\hat r_{\hat X}(\hat x)-1}(\hat x))=e^{-\hat\phi(\hat f^j(\hat x))} d\hat m_{\phi}(\hat f^j(\hat x)),
\end{align*}
as required.

For the proof of (b), %let $\hat X_i$ be some cylinder of $(\hat X,R_{\hat X})$.
 for $\hat x\in \hat X$, by definition $d\hat m_\phi'(R_{\hat X}(\hat x))=e^{-\hat\Phi(\hat x)}d\hat m_\phi'(\hat x)$.  Let $\hat X'$ be some domain in $\hat X$ contained in some single domain $D\in \D$ (this is not a necessary step if the inducing scheme is of type A). This implies that $m_\phi':=\hat m_\phi'\circ\pi_{\hat X'}^{-1}$ is $\Phi$-conformal after rescaling.  As in Proposition~\ref{prop:eqgen}, there is only one $\Phi$-conformal measure for $(X,F)$, which implies that $\hat m_\phi'=\hat m_\phi$ up to a rescaling.
\end{proof}

Given $\hat X\subset \hat I$, we consider the system $(\hat X,R_{\hat X})$ where $R_{\hat X}$ is the first return map to $\hat X$.  The measure $\hat\mu_\phi$ is an invariant measure for $(\hat X,R_{\hat X})$, see \cite{Kbook}.  Adding Kac's Lemma to \eqref{eq:lift}, for any $\hat A\subset \hat I$ we have \begin{equation}\hat\mu_\phi(\hat A):=\sum_i\sum_{0\le k\le r_{\hat X}|_{\hat X_i}-1}\hat\mu_{\phi}(\hat f^{-k}(\hat A)\cap \hat X_i).\label{eq:saturate}
\end{equation}
This means we can compare $\hat m_\phi$ and $\hat\mu_\phi$ on domains $\hat f^j(\hat X_i)$, for $0\le k\le r_{\hat X}|_{\hat X_i}-1$, in a relatively simple way.

We will project the measure $\hat m_\phi$ to $I$.  Although it is possible to show that for many potentials we consider, $\hat m_\phi(\hat I)<\infty$, we allow the possibility that our conformal measures are infinite.  This leaves the possibility to extend this theory to a wider class of measures open.  So in the following lemma, we use another way to project $\hat m_\phi$.

\begin{lemma}
Suppose that $\hat Y\subset \hat I_{\T}$ is so that $\hat Y=\sqcup_n\hat Y_n$ for $Y_n$ an interval contained in a single domain $D_{Y_n}\in \D_{\T}$ and $\pi:\hat Y\to I$ is a bijection.  Then for $\nu_\phi:=\hat m_\phi\circ\pi|_{\hat Y}^{-1}$, we have $\nu_\phi(I)<\infty$.  Moreover, $m_\phi:=\frac{\nu_\phi}{\nu(I)}$ is a conformal measure for $(I,f,\phi)$, and $m_\phi$ is independent of $\hat Y$.
\label{lem:conf proj}
\end{lemma}

\begin{proof}  We first prove that $\nu_\phi$ is independent of $\hat Y$, up to rescaling.  In doing so, the $\phi$-sigma-conformal property of $\nu_\phi$ become clear.  The we show that $\nu_\phi(I)<\infty$.

Let us pick some $\hat Y$, and let $\nu_\phi$ be as in the statement of the lemma.
Let $x\notin\cup_{n\in \N}f^n(\crit)$.  Suppose that $\hat x_1, \hat x_2$ have $\pi(\hat x_1)=\pi(\hat x_2)=x$.    By our condition on $x$, we have $\hat x_i\notin \bd\hat I$ for $i=1,2$.  We denote $D_1,D_2\in \D$ to be the domains containing $x_1,x_2$ respectively.  The independence of the measure from $\hat Y$ follows if we can show for any neighbourhood $U$ of $x$ such that for $\hat U_i:=\pi^{-1}(U)\cap D_i$ such that $\hat U_i\Subset D_i$ for $i=1,2$, we have $\hat m_\phi(\hat U_1)=\hat m_\phi(\hat U_2)$.

As in \cite{Kellift} there exists $n\ge 0$ so that $\hat f^n(\hat x_1)=\hat f^n(\hat x_2)$.  Since we are only interested in the infinitesimal properties of our measures, we may assume that the same is true of $\hat U_1$ and $\hat U_2$, \ie $\hat f^n(\hat U_1)=\hat f^n(\hat U_2)$.  Therefore $\hat m_\phi(\hat f^n(\hat U_1))=\int_{\hat U_1} e^{-\hat\phi_n}~d\hat m_\phi$.  Since $\hat m_\phi(\hat f^n(\hat U_1))=\hat m_\phi(\hat f^n(\hat U_2))$ and $\hat \phi=\phi\circ\pi$, we have $\hat m_\phi(\hat U_1)=\hat m_\phi(\hat U_2)$, as required.  So it only remains to show $\nu_\phi(I)<\infty$.

By the above, the $\hat\phi$-sigma-conformality of $\hat m_\phi$ passes to $\phi$-sigma-conformality of $\nu_\phi$.  We can pick $U\subset I$ such that $U=\pi(\hat U)$ for some $\hat U \subset D\in \D_{\T}$.  Recall that $m_\phi$ was obtained from a conformal measure $m_\Phi$ for some inducing scheme $(X,F)$.  We may assume that $\hat U$ is such that $\hat U\subset \hat f^k(\hat X_i)\cap D$ for some $0\le k\le r_{\hat X}|_{\hat X_i}-1$ and some $D\in \D$.  This implies that $\hat m_\phi(\hat U)<\infty$, and so $\nu_\phi(U)<\infty$.  Since $f$ is in $\F$, it is locally eventually onto, \ie for any small open interval $W\subset I$ there exists $n\in \N$ so that $f^n(W) \supset \Omega$.  Therefore there exists $n$ so that $f^n(U)\supset I$.  Then by the $\phi$-sigma-conformality of $\nu_\phi$, we have
$$\nu_\phi(I)=\nu_\phi(f^n(U))=\int_Ue^{-\phi_n}~d\nu_\phi \le \nu_\phi(U)e^{-\inf\phi_n}<\infty.$$ Hence $m_\phi$ is conformal.
\end{proof}

Note that combining Lemmas~\ref{lem:sig conf} and \ref{lem:conf proj}, we deduce that $m_\phi$ is independent of the inducing scheme that produced it.  We next consider the density.

\begin{lemma}
For $\phi\in SVI$ satisfying \eqref{eq:range} and $P(\phi)=0$, $\frac{d\mu_\phi}{dm_\phi}$ is uniformly bounded above.
\label{lem:dens bound}
\end{lemma}

\begin{proof}   Suppose that $\frac{d\mu_\phi}{dm_\phi}(x)>0$.  We let $\pi^{-1}(x)=\{\hat x_1,\hat x_2, \ldots\}$, where the ordering is by the level, \ie $\level(\hat x_{j+1})\ge \level(\hat x_j)$ for all $j\in \N$.  Then since $\mu_\phi=\hat\mu_\phi\circ\pi^{-1},$ $$
\frac{d\mu_\phi}{dm_\phi}(x) = \sum_{j=1}^\infty\frac{d\hat\mu_\phi}{dm_\phi\circ\pi}(\hat x_j).$$
We will use this fact allied to equation \eqref{eq:saturate} for return maps on the Hofbauer extension, and the bounded distortion of the measures for these first return maps to get the bound on the density. We note that since for any $R\in \N$, there are at most $2\#\crit$ domains of $\D$ of level $R$ (see for example \cite[Chapter 9]{BrBr}), there can be at most $2\#\crit$ elements $\hat x_j$ of the same level.

We let $(X,F)$ be a type A inducing scheme with induced potential $\Phi:X \to \R$.  Let $\hat X$ be the interval in $\hat I$ for which the first return map $R_{\hat X}$ defines the inducing scheme $(X,F)$.
Recall that $\mu_\Phi$ can be represented as $\frac{\hat\mu_\phi\circ\pi|_{\hat X}^{-1}}{\hat\mu_\phi(\hat X)}$ and by Lemma~\ref{lem:conf proj}, we can express $m_\Phi$ as $\frac{m_\phi}{m_\phi(X)}$.  Moreover as in Proposition~\ref{prop:eqgen} there exists $C_\Phi>0$ so that $\frac{d\mu_\Phi}{dm_\Phi}\le C_\Phi$.

Since $R_{\hat X}$ is a first return map, for each $i$ there exists at most one point $\hat x_{j,i}$ in $\hat X_i$ so that $\hat f^k(\hat x_{j,i})=\hat x_j$ for $0\le k< r_{\hat X}|_{\hat X_i}$.  We denote this value $k$ by $r_{j,i}$.  Let $k_j:=\inf\{r_{j,i}:i\in \N\}$.

By \eqref{eq:saturate}, $d\hat\mu_\phi(\hat x_j) = \sum_id\hat\mu_\phi(\hat x_{j,i})$.
By conformality, for each $i$, $$d\hat m_\phi(\hat x_j)= e^{-\hat\phi_{r_{j,i}}(\hat x_{j,i})}~d\hat m_\phi(\hat x_{j,i}) \ge e^{-\sup\phi_{r_{j,i}}}~d\hat m_\phi(\hat x_{j,i}).$$
Therefore, letting $x_{j,i}=\pi(\hat x_{j,i})$,
\begin{align*}
\frac{d\hat\mu_\phi}{d\hat m_\phi}(\hat x_j) & \le \sum_i\frac{d\hat\mu_\phi}{d\hat m_\phi}(\hat x_{j,i}) e^{\sup\phi_{r_{j,i}}} \le \left(\frac{m_\phi(X)}{\hat\mu_\phi(\hat X)}\right)\sum_i\frac{d\mu_\Phi}{ dm_\Phi}(x_{j,i}) e^{\sup\phi_{r_{j,i}}}\\
&\le C_\Phi\left(\frac{m_\phi(X)}{\hat\mu_\phi(\hat X)}\right)\sum_i e^{\sup\phi_{r_{j,i}}} \le C_\Phi\left(\frac{m_\phi(X)}{\hat\mu_\phi(\hat X)}\right)\sum_n \#\{i:r_{j,i}=n\} e^{n\sup\phi}.
\end{align*}

By \cite{Htop}, if $\level(\hat x_j)=R$ then there exist $C>0$ and $\gamma(R)>0$ so that $\gamma(R)\to 0$ as $R\to \infty$ and the number of $n$-paths terminating at $D_{\hat x_j}\in \D$ at most $Ce^{n\gamma(R)}$.  Then $\#\{i:r_{j,i}=n\} \le Ce^{n\gamma(\level(\hat x_j))}$.  Also $k_j\ge \level(\hat x_j)-\level(\hat X)$. Therefore, \begin{align*}\frac{d\hat\mu_\phi}{d\hat m_\phi}(\hat x_j) & \le CC_\Phi \left(\frac{m_\phi(X)}{\hat\mu_\phi(\hat X)}\right) \sum_{n\ge k_j} e^{n\left(\gamma(\level(\hat x_j))+\sup\phi\right)} \\
&\le CC_\Phi \left(\frac{m_\phi(X)}{\hat\mu_\phi(\hat X)}\right) e^{\left(\level(\hat x_j)-\level(\hat X)\right)\left(\gamma(\level(\hat x_j))+\sup\phi\right)}\sum_{n\ge 0} e^{n\left(\gamma(\level(\hat x_j))+\sup\phi\right)}.\end{align*}
Since, as in Lemma~\ref{lem:press 0}, our conditions on $\phi$ ensure that $\sup\phi<0$, there exists $\kappa>0$, and $j_0\in \N$ so that $\gamma(\level(\hat x_j))+\sup\phi<-\kappa$ for all $j\ge j_0$.    Since there are at most $2\#\crit$ points $\hat x_j$ of any given level $R$, there are only finitely many $j$ with $\level(\hat x_j)-\level(\hat X)\le 0$.  Moreover, there exists $C'>0$ so that
$$\frac{d\mu_\phi}{dm_\phi}(x) \le \sum_{j=1}^{j_0-1}\frac{d\hat\mu_\phi}{dm_\phi\circ\pi}(\hat x_j) +\sum_{j=j_0}^\infty \frac{d\hat\mu_\phi}{dm_\phi\circ\pi}(\hat x_j) \le C'+C'\sum_{j=j_0}^\infty e^{-j\kappa}$$
which is uniformly bounded.
\end{proof}

\begin{proof}[Proof of Proposition~\ref{prop:conf}']
The existence of the conformal measure $m_\phi$ is proved in the above lemmas.  Lemma~\ref{lem:dens bound} implies that the density $\frac{d\mu_\phi}{dm_\phi}$ is uniformly bounded above. The lower bound follows by a standard argument, which we give for completeness. Proposition~\ref{prop:eqgen} implies that we can take a type A inducing scheme $(X, F,\Phi)$ so that $\frac{d\mu_\Phi}{dm_\Phi}$ is uniformly bounded below by some $C_\Phi^{-1}\in (0,\infty)$.  Also, Lemma~\ref{lem:sig conf} implies that $\frac{m_\phi}{m_\phi(X)}=m_\Phi$.  Since, as in the proof of Lemma~\ref{lem:conf proj}, $(I,f)$ is locally eventually onto, there exists $n\in\N$ so that $f^n(X)\subset \Omega$.  So for a small interval $A\subset \Omega$,  there exists some $A_i\subset X_i$ so that $f^k(A_i)=A$ for some $0\le k\le n$.  Then \eqref{eq:lift} implies that
$$\frac{\mu_\phi(A)}{m_\phi(A)}\ge \frac{\mu_\phi(A_i)}{m_\phi(A_i)}e^{\inf\phi_n} \ge \left(\frac{m_\phi(X)}{\int\tau~d\mu_\Phi}\right)\hspace{-1mm} \left(\frac{\mu_\Phi(A_i)}{m_\Phi(A_i)}\right)\hspace{-1mm} e^{\inf\phi_n}\ge
\left(\frac{m_\phi(X)}{\int\tau~d\mu_\Phi}\right)\hspace{-1mm} \left(\frac{e^{\inf\phi_n}}{C_\Phi}\right).$$ Hence $\frac{d\mu_\phi}{dm_\phi}$ is uniformly bounded below.
\end{proof}

\begin{lemma}
Suppose that $f\in \F$ satisfies \eqref{eq:grow} and $\phi\in SVI$.  Then there exists $\eps>0$ so that for any inducing scheme $(X,F)\in SCover^B(\eps)$, the induced potential $\Phi$ has $P(\Phi)=0$.
\label{lem:press 0}
\end{lemma}

\begin{proof}
We will apply Case 3 of \cite[Proposition 1]{BTeqnat}.  Firstly we need to show that $Z_0(\Phi)<\infty$.  By   Proposition~\ref{prop:conf}' there exists a conformal measure $m_\phi$, coming from an inducing scheme of type A in Proposition~\ref{prop:eqgen}'.
By the $\phi$-conformality of $m_\phi$ and the local H\"older continuity of $\Phi$, as in Proposition~\ref{prop:eqgen}(b), there exists $C>0$ so that $Z_0^*(\Phi) \le C \sum_i\tau_i m_\phi(X_i).$
Then by Proposition~\ref{prop:conf}' and the facts that $(X,F)$ was generated by a first return map to some $\hat X$ and $\mu_\phi=\hat\mu_\phi\circ\pi^{-1}$,
$$Z_0^*(\Phi) \le CC_\phi' \sum_i\tau_i \mu_\phi(X_i)=CC_\phi' \sum_ir_{\hat X}|_{\hat X_i} \hat\mu_\phi(\hat X_i).$$
By Kac's Lemma this is bounded.

Now the fact that $\mu_\phi$ is compatible to $(X,F)$ follows simply, see for example Claim 1 in the proof of \cite[Proposition 2]{BTeqnat}.  Then Case 3 of \cite[Proposition 1]{BTeqnat} implies $P(\Phi)=0$.
\end{proof}

\begin{proof}[Proof of Proposition~\ref{prop:induced density}']
Suppose that $(X,F)$ is an inducing scheme as in the statement, with induced potential $\Phi$.  If $(X,F)$ is of type A then by Lemma~\ref{lem:sig conf}, the measure $m_\phi$ works as a conformal measure for $(X,F,\Phi)$, up to renormalisation.  By Proposition~\ref{prop:eqgen}(c), $m_\phi$ is in fact equal to $m_\Phi$ up to renormalisation.  By Lemma~\ref{lem:press 0}, this is also true for type B inducing schemes.  Since by Proposition~\ref{prop:conf}', $\frac{d\mu_\phi}{dm_\phi}$ is bounded above and below, and as in Proposition~\ref{prop:eqgen}, we have $\frac1{C_\Phi}\le \frac{d\mu_\Phi}{dm_\Phi}\le C_\Phi$, this implies that $\frac{d\mu_\Phi}{d\mu_\phi}$ is also uniformly bounded above and below.
\end{proof}

\end{document}